\documentclass[11pt]{article}
\usepackage{amsmath, amsthm, amssymb, bbm, setspace,bigints, amsfonts, hyperref}
\usepackage[round,authoryear]{natbib}
\usepackage[margin=1 in]{geometry}
\usepackage{caption}
\usepackage{subcaption}
\usepackage[toc,page]{appendix}
\usepackage[pdftex]{graphicx}
\usepackage{pdfpages}
\usepackage{epstopdf}
\usepackage{booktabs}
\usepackage{comment}
\usepackage{float}
\usepackage{mathtools}
\usepackage{booktabs,longtable}
\mathtoolsset{showonlyrefs}
\usepackage[ruled,vlined]{algorithm2e}
\usepackage{authblk}

\newcommand{\be} {\begin{eqnarray*}}
\newcommand{\ee} {\end{eqnarray*}}

\newcommand{\diag} {\mathrm{diag}}

\newcommand{\bk} {{\bf k}}
\newcommand{\bi} {{\bf i}}
\newcommand{\bj} {{\bf j}}
\newcommand{\boldm} {{\bf m}}
\newcommand{\bell} {{\bf l}}
\newcommand{\bt} {{\bf t}}
\newcommand{\bs} {{\bf s}}

\newcommand{\bx} {{\bf x}}

\newcommand{\cN} {{\mathcal{N}}}

\newcommand \bbB{\mathbb{B}}

\newcommand \bbE{\mathbb{E }}
\newcommand \bbH{\mathbb{H}}

\newcommand \bbN{\mathbb{N}}
\newcommand \bbR{\mathbb{R}}

\newcommand \cD{\mathcal{D}}
\newcommand \cE{\mathcal{E}}
\newcommand \cF{\mathcal{F}}

\newcommand \cL{\mathcal{L}}
\newcommand \cI{\mathcal{I}}

\newcommand \cR{\mathcal{R}}

\newcommand{\cT}{\mathcal{T}}

\newcommand{\id}{\mathrm{id}}

\newcommand{\cB}{\mathcal{B}}
\newcommand{\cG}{\mathcal{G}}
\newcommand{\cJ}{\mathcal{J}}

\newcommand \TS{\Tilde{S}}
\newcommand \Besova{B^{\alpha}_{\infty, \infty}}

\newcommand \GP{\mathrm{GP}}

\newcommand \Data{{\mathcal{D}_n}}

\newcommand{\ind}{\mathrm{1}}

\newtheorem{theorem}{Theorem}[section]
\newtheorem{lemma}[theorem]{Lemma}

\newtheorem{corollary}[theorem]{Corollary}

\newcommand{\papertitle}{Adaptive Resolution for Finite-Rank Gaussian Processes}

\title{\papertitle}

\author[1]{Jaehoan Kim\thanks{jaehoan.kim@duke.edu}}
\author[3]{Anirban Bhattacharya \thanks{anirbanb@stat.tamu.edu}}
\author[2]{Debdeep Pati\thanks{dpati2@wisc.edu}}
\affil[1]{Department of Statistical Science, Duke University}
\affil[2]{Department of Statistics, University of Wisconsin-Madison}
\affil[3]{Department of Statistics, Texas A\&M University}

\begin{document}
\maketitle

\begin{abstract}
Finite-rank approximations are widely used to scale Gaussian process (GP) regression, but their posterior behavior can differ from that of the corresponding parent GP prior. We study a class of finite-rank GP priors built from locally supported basis expansions with dependent Gaussian coefficients. Our framework covers finite-element approximations based on the stochastic partial differential equation (SPDE) representation of Mat\'ern GPs and regular-grid GP interpolation schemes. We show that, with a suitable prior on the resolution parameter $N$, these finite-rank expansions inherit the same posterior contraction rate as the corresponding parent GP prior under the same bandwidth specification used for that parent prior. Consequently, the interpolation construction under a squared-exponential parent GP attains the minimax-optimal rate up to logarithmic factors under a hierarchical prior on the bandwidth parameter and on $N$, while the SPDE construction attains the same rate under a bandwidth scaling depending on the sample size and the smoothness of the true function, together with a prior on $N$. We also develop a posterior sampler for the hierarchical interpolation model that jointly updates the resolution and bandwidth parameters, and we provide numerical studies that support the theory.
\end{abstract}


\section{Introduction}
Gaussian process (GP) methods are a standard tool for Bayesian learning of function-valued parameters, but exact implementation requires $O(n^2)$ memory and $O(n^3)$ floating point operations. This quickly becomes infeasible even when $n$ is only moderately large. As a result, scalable GP approximations have received substantial attention over the last few decades. Broadly speaking, these methods fall into three main categories: sparse approximations, low-rank approximations, and distributed computing methods. For methods that combine features of the first two categories, see \citet{sun2012geostatistics, heaton2019case} for reviews.

Despite this broad methodological development, there is still limited theoretical understanding of when such approximations retain the statistical guarantees of the corresponding parent GP prior. Many popular sparse and low-rank methods can be viewed through the lens of Vecchia's approximation \citep{vecchia1988estimation}, which replaces the full Gaussian likelihood by a product of lower-dimensional conditional distributions. Other examples of sparse approximations include pseudo-likelihood methods \citep{varin2011overview}, spectral methods \citep{fuentes2007approximate}, and covariance tapering \citep{furrer2006covariance}. Representative low-rank constructions include process convolutions \citep{higdon2002space}, fixed rank kriging \citep{cressie2008fixed}, predictive processes \citep{banerjee2008gaussian}, lattice kriging \citep{nychka2015multiresolution}, and the stochastic partial differential equation (SPDE) approach \citep{lindgren2011explicit}. Distributed approaches aggregate local experts trained on subsets of the data; see, for example, product-of-experts \citep{deisenroth2015distributed, cao2014generalized} and hierarchical mixture-of-experts \citep{ng2014hierarchical}.

Recent theory has started to address this question for several scalable GP procedures. \citet{stein2004approximating} showed that the Vecchia likelihood is the best Kullback--Leibler approximation of the target Gaussian distribution under a sparsity constraint on the Cholesky factor of the precision matrix. \citet{kang2024asymptotic} established that Vecchia prediction can approximate exact GP prediction when the conditioning set size grows suitably. Posterior contraction theory has also been developed for Vecchia Gaussian processes and related directed-graph constructions \citep{szabo2024vecchia, zhu2024radial}. For distributed Bayesian procedures, posterior convergence and adaptation have been studied in \citep{szabo2020adaptive, szabo2019asymptotic, guhaniyogi2023distributed}. For sparse variational approximations, \citet{nieman2022contraction} showed that variational posteriors can attain minimax-optimal contraction rates when the number of inducing variables is chosen appropriately; see also \citet{nieman2025adaptive} for recent results on adaptive hyperparameter selection through variational Bayes.

A common feature of these results is that statistical optimality depends critically on hyperparameter choice. Examples include truncation levels, smoothness or range parameters in the covariance kernel, and approximation-specific tuning parameters such as the number of inducing variables or the size of local conditioning sets. For finite-rank approximations, a central approximation-specific parameter is the resolution $N$, which determines the dimension of the expansion and hence the support of the induced prior. This motivates studying finite-rank constructions in a way that keeps the role of $N$ explicit.

In this article, we study a broad class of finite-rank GP approximations of the form $f_N(x) = \sum_{j \in \cJ_N} w_j \psi_j(x)$, where $\cJ_N$ is an index set and $\{\psi_j : j \in \cJ_N\}$ is a collection of basis functions on $[0,1]^d$. We endow the coefficient vector $w^N$ with a dependent Gaussian prior of the form $w^N \mid N,\theta \sim N(0,\Sigma_{N,\theta})$, where $\Sigma_{N,\theta}$ is a non-diagonal positive definite covariance matrix and $\theta$ denotes the kernel hyperparameters. The covariance $\Sigma_{N,\theta}$ is induced from the covariance kernel $K_\theta$ of an underlying full Gaussian process $f \sim \mathrm{GP}(0, K_\theta)$, which we call the \emph{parent GP}. The way $\Sigma_{N,\theta}$ is obtained from $K_\theta$ depends on the construction. To complete the prior specification, we later place a prior on $N$,
while $\theta$ is either taken to be a sample-size-dependent sequence inherited from the corresponding parent-GP theory or endowed with an independent hyperprior. This
defines the general class of finite-rank GP priors considered throughout the paper;
the formal specification is given in Section~\ref{sec:2methods}.

A prominent example of this framework is process convolution \citep{higdon2002space}, where $\psi_j$ are kernel basis functions and the coefficients $w_j$ are independent Gaussian. 
With locally supported bases and independent coefficients, however, the global smoothness of $f_N$ is limited by the fixed regularity of the basis functions. As a result, sufficiently smooth targets may not be well approximated \citep{10.1214/12-EJS735}. A natural alternative is to endow the coefficient vector $w^N$ with a dependent Gaussian prior. Then $f_N$ can inherit the global regularity of its parent GP while retaining the computational advantages induced by local support.

We study two representative constructions within the form of \eqref{eqn:f_approximation}--\eqref{eqn:w_prior}. The first is the finite-element SPDE approximation \citep{lindgren2011explicit}, where the coefficient vector is Gaussian with sparse precision structure induced by the SPDE representation. Related finite-element representations of Gaussian processes have also been studied from the perspective of balancing numerical and statistical accuracy \citep{sanz2022finite}. The second is a lattice-based GP interpolation approach used in computer emulation \citep{maatouk2017gaussian, zhou2019reexamining, zhou2024mass, maatouk2025bayesian}, where the coefficients on a regular grid are assigned a Gaussian prior induced from a parent GP, typically together with a bandwidth parameter.

A key point is that the role of $N$ is not merely computational. Conditional on fixed $(N,\theta)$, the prior is supported on $\mathrm{span}\{\psi_j : j \in \cJ_N\}$. Hence for small $\epsilon > 0$, the event $\{\|f_N - f^*\|_\infty < \epsilon\}$ can have negligible, or even zero, conditional prior probability unless $N$ is sufficiently large. In general, the required size of $N$ depends on additional hyperparameters through $\theta$. Thus, even when the parent GP prior is statistically well understood, a finite-rank approximation need not inherit that behavior unless the prior on $N$ places enough mass on sufficiently large resolutions.

Our main result is that a suitable prior on the resolution $N$ allows the finite-rank approximation to inherit the posterior contraction behavior of the corresponding parent GP prior. We call this an {\it adaptive-resolution result}: the approximation resolution is learned through the prior on $N$, while the bandwidth or range parameter is treated according to the parent-GP theory used for the corresponding construction. For the interpolation construction with a squared-exponential parent GP, we place hierarchical priors on both the bandwidth parameter and $N$, yielding the minimax-optimal posterior contraction rate up to logarithmic factors. For the SPDE construction, we place a prior on $N$ and use the bandwidth scaling inherited from the parent Mat\'ern contraction theory. Thus the common contribution across the two constructions is adaptation over the finite-rank resolution, showing that the finite-dimensional approximation can preserve the parent-GP contraction rate when sufficiently large resolutions receive enough prior mass.

The SPDE result also connects two related strands of literature. On the one hand, \citet{sanz2022finite} study finite-element GP representations from the perspective of numerical and statistical approximation error at fixed resolution. On the other hand, \citet{fang2025posterior} establish posterior contraction results for rescaled and hierarchical Mat\'ern GP priors. Our SPDE theorem links these viewpoints by showing that, under the bandwidth scaling used in the parent Mat\'ern theory, a suitable prior on the finite-element resolution preserves the parent contraction rate for the finite-rank approximation.

While the adaptive-resolution hierarchy is theoretically attractive, it, however, introduces a computational challenge because the coefficient vector $w^N$ changes dimension with $N$. We address this by integrating out $w^N$ when updating $N$ and the kernel hyperparameters, and then sampling $w^N$ conditionally from its Gaussian posterior. This avoids dimension-changing reversible-jump updates \citep{green1995reversible}. The resulting computations exploit the local support of the basis functions through sparse design matrices. Thus, the data-dependent part of each iteration scales linearly in sample size $n$, and the remaining cost is governed by the coefficient dimension $|\cJ_N|$. This involves $O(|\cJ_N|^3)$ calculation for the GP interpolation method, and is further reduced for the SPDE method due to sparsity in the posterior precision matrix.

The rest of the article is organized as follows. In Section \ref{sec:2methods}, we introduce the two approximation frameworks based on the SPDE formulation and lattice-based interpolation. Section \ref{sec:RKHS} outlines the reproducing kernel Hilbert space structure associated with \eqref{eqn:f_approximation}. Posterior contraction theory for the SPDE-based approximation and the interpolation-based approximation is developed in Sections \ref{ssec:SPDEt} and \ref{ssec:GPIt}, respectively. Section \ref{sec:comp} discusses computational strategies, followed by numerical experiments in Section \ref{sec:sims}. We conclude with a discussion in Section \ref{sec:disc}. Proofs of theoretical results are provided in a supplemental document. 
\section{Finite rank approximation of Gaussian processes}\label{sec:2methods}
We focus on the nonparametric regression model with Gaussian errors
\begin{align}\label{eqn: fixed regression model}
     y_i = f(x_i) + \epsilon_i, \quad \epsilon_i \overset{ind.}\sim N(0, \sigma^2),
\end{align}
with the observed data $\cD_n = \{(x_1, y_1), \ldots, (x_n, y_n)\}$. Throughout the article, for simplicity of exposition, we assume $\sigma$ to be known and the design points $\{x_i\}_{i=1}^n$ to be fixed. The infinite-dimensional parameter of interest is $f: [0,1]^d \to \mathbb{R}$.
 
Throughout, parent Gaussian process priors are written as $f \sim \mbox{GP}(0, K)$, where $K$ denotes the covariance kernel. We represent the regression function $f$ in \eqref{eqn: fixed regression model} by a finite basis expansion form
\begin{equation}\label{eqn:f_approximation}
    f_N(x) = \sum_{\bj\in \cJ_N} w_\bj \psi_\bj(x),
\end{equation}
where $\{\psi_\bj: \bj\in \cJ_N\}$ is a collection of basis functions on $[0,1]^d$ with local support and $\cJ_N := \{(j_1,\ldots,j_d): j_k\in\{0,1,\ldots,N\}\}$.
We endow the coefficient vector $w^N=(w_\bj)_{\bj\in\cJ_N}$ with a dependent Gaussian prior of the form
\begin{equation}\label{eqn:w_prior}
    w^N \mid N,\theta \sim N(0,\Sigma_{N,\theta}),
\end{equation}
where $\Sigma_{N,\theta}$ is a non-diagonal positive definite covariance matrix induced from a parent GP, and $\theta$ denotes the hyperparameters involved in the covariance kernel of the parent GP. In both constructions studied below, $\theta$ is an inverse bandwidth parameter $\kappa$. The treatment of $\kappa$ differs across the two methods: in the SPDE approach, $\kappa=\kappa_n$ with $\kappa_n$ depending on the sample size and smoothness of the true function, inherited from the parent Mat\'ern GP theory, while in the GP interpolation approach, $\kappa$ is assigned a hyperprior. In both cases, a prior on $N$ is placed for the adaptive resolution. We refer to a prior of the form \eqref{eqn:f_approximation}--\eqref{eqn:w_prior} with locally supported $\{\psi_\bj\}$ and non-diagonal $\Sigma_{N,\theta}$ as a {\it finite rank GP with dependent coefficients and locally supported bases} (\textbf{d-frGPloc}).

\subsection{Finite element approximation to GP as a solution to SPDE}\label{ssec:SPDE}
 
As our first candidate for d-frGPloc, we consider the SPDE-based approach of \citet{lindgren2011explicit}. Recall that the stationary Mat\'{e}rn covariance kernel with smoothness parameter $\nu$, marginal variance $\tau^2$, and inverse length-scale parameter $\kappa$ is
\begin{align}\label{eqn: matern cov kernel}
K(x_1, x_2)
= \tau^2 \frac{2^{1-\nu}}{\Gamma(\nu)} \, (\kappa \|x_1 - x_2\|_2)^\nu \,
K_{\nu}(\kappa \|x_1 - x_2\|_2),
\end{align}
where $K_\nu$ is the modified Bessel function of the second kind. A well-known result from \citet{whittle1963stochastic} is that on $\bbR^d$, a Mat\'ern GP can be characterized as the stationary solution to an SPDE of the form
\begin{align}\label{eqn: SPDE for Matern}
    (\kappa^2 - \Delta)^{\beta/2} f = \kappa^{\beta - d/2} W,
\end{align}
where $\beta=\nu+d/2$, $W(\cdot)$ is Gaussian white noise, and $\Delta$ denotes the Laplacian operator $\Delta(f)=\sum_{i=1}^d \frac{\partial^2 f}{\partial x_i^2}$. The marginal variance parameter $\tau^2$ in \eqref{eqn: matern cov kernel} acts only as a scaling factor and will be suppressed in the SPDE representation below.

Starting from this result, the SPDE-based finite-rank prior is constructed as follows. We first use the SPDE representation to define the parent GP, and then approximate its weak solution in a finite-dimensional space spanned by locally supported finite-element basis functions. The resulting approximation has the form of the general expansion in equation \eqref{eqn:f_approximation}, with a sparse precision matrix for the coefficient vector. See \citet{lindgren2011explicit, sanz2022finite} for details on the finite-element SPDE construction.
 
In this article, we take $\beta\in\bbN$ and consider the bounded domain $D=[0,1]^d$. For an SPDE constructed on a bounded domain $D \subset \bbR^d$, the boundary condition is required for its solution to be well-defined. In this article, we impose homogeneous Neumann boundary conditions, as in Section~2.3 of \citet{lindgren2011explicit} and Section~3 of \citet{sanz2022finite}.
Other boundary conditions may also be considered, but they lead to slightly different finite-element solutions and are not pursued here. Under this choice, the operator $\cL_\kappa=\kappa^2\,\mbox{id}-\Delta$ is self-adjoint and positive on $L^2(D)$, and $(\kappa^2-\Delta)^{\beta/2}$ is understood via spectral powers of $\cL_\kappa$. When $\beta$ is even, this coincides with $\beta/2$ iterated applications of $\cL_\kappa$.
 
For the given SPDE, a random function $g$ is called a generalized solution if it satisfies
\begin{align}\label{eqn: weak solution condition}
   \langle \phi, (\kappa^2 - \Delta)^{\beta/2} g \rangle \overset{d}{=} \langle \phi, \kappa^{\beta - d/2} W \rangle
\end{align}
for every smooth test function $\phi$ and $\langle \cdot, \cdot \rangle$ denoting inner product on $D$. Following the Galerkin construction in section 2.3 of\citet{lindgren2011explicit}, we approximate the solution by considering a finite collection of test functions $\{\phi_\bj\}_{\bj\in\cJ_N}$:
\begin{align}\label{eqn: galerkin condition}
   \langle \phi_\bj, (\kappa^2 - \Delta)^{\beta/2} g \rangle \overset{d}{=} \langle \phi_\bj, \kappa^{\beta - d/2} W \rangle,
   \quad \bj\in\cJ_N.
\end{align}
If one sets $\phi_\bj=\psi_\bj$ for $\beta \ge 2$ and $\phi_\bj = (\kappa^2 - \Delta)^{1/2} \psi_\bj$ for $\beta = 1$, then the resulting Galerkin solution takes the form \eqref{eqn:f_approximation}.
We then define the consistent finite-element mass and stiffness matrices by
\begin{equation}
C^{(c)}_{\bi,\bj}=\langle \psi_{\bi},\psi_{\bj}\rangle,
\qquad
G^{(c)}_{\bi,\bj}= \langle \nabla \psi_{\bi}, \nabla \psi_{\bj}\rangle,
\qquad \bi,\bj\in\cJ_N .
\label{eqn: CG defn}
\end{equation}

Let $C=\diag(C^{(c)}\ind)$ be the lumped mass matrix, where $\ind$ is a vector of ones of appropriate dimension. We write $G$ for the corresponding mass-lumped stiffness matrix. We define
\begin{equation*}
L_N=\kappa^2 I+C^{-1}G.
\end{equation*}
Here, $L_N$ is a mass-lumped finite-element coefficient-coordinate discretization of $\cL_\kappa = \kappa^2\id-\Delta$. The precise transition from the consistent matrices to the mass-lumped tensor-product matrices is given in Appendix \ref{appendix:SPDE_technical}.

From the local support of the finite-element basis, the resulting $L_N$ becomes a sparse matrix. Note that the sparsity is obtained by replacing $C^{(c)}$ with the lumped matrix $C$, since although $C^{(c)}$ is a sparse matrix, $\big( C^{(c)}\big)^{-1}$ may be a dense matrix.

For $\beta\in\bbN$, define the scale-corrected precision matrix
\begin{equation}
Q=
\kappa^{-(2\beta-d)}C \, L_N^\beta.
\label{eqn:Q_def_main}
\end{equation}
Then the approximate finite-rank solution of the SPDE in equation \eqref{eqn: SPDE for Matern} is written in the form of equation \eqref{eqn:f_approximation}, with the prior in equation \eqref{eqn:w_prior} where $\theta=\kappa$ and
\begin{equation*}
\Sigma_{N,\theta}=Q_{N, \theta}^{-1},
\end{equation*}
with $Q_{N, \theta} = Q$.
We henceforth refer to inference based on this prior as the {\it SPDE method}. The local support of the finite-element basis makes $L_N$ and $Q$ sparse for integer $\beta$. The sparsity in $Q$ is the main computational advantage of the SPDE construction.

As an illustration, consider $D=[0,1]$ and the piecewise linear hat basis $\{\psi_j\}$,
\begin{align}\label{eqn: definition of linear basis function}
    \psi_j(x) = \Big(1 - N\Big|x - \frac{j}{N}\Big|\Big)^+,
    \quad j=0,\ldots,N,
\end{align}
where $x^+ := \max(x, 0)$. The lumped mass matrix $C$ and stiffness matrix $G$ then assume the form
\begin{align}\label{eq:CG}
    C = \frac{1}{2N}\begin{bmatrix}
    1 & 0 & \cdots & 0 & 0\\
    0 & 2 & \cdots & 0 & 0 \\
    \vdots & \vdots & \ddots & \vdots & \vdots \\
    0 & 0 & \cdots & 2 & 0 \\
    0 & 0 & \cdots & 0 & 1
    \end{bmatrix},\quad
    G = N\begin{bmatrix}
    1 & -1 & 0 & \cdots & 0 & 0 & 0\\
    -1 & 2 & -1 & \cdots & 0 & 0 & 0\\
    \vdots & \vdots & \vdots &  \ddots & \vdots & \vdots & \vdots \\
    0 & 0 & 0 & \cdots & -1 & 2 & -1 \\
    0 & 0 & 0 & \cdots & 0 & -1 & 1
    \end{bmatrix}.
\end{align}
The bandwidth of $Q$ is determined by $\beta$.

\subsection{Interpolating GP on a regular grid}\label{ssec:GPI}
Our second construction starts from a parent GP prior and approximates it by interpolation on a regular grid \cite{maatouk2017gaussian,zhou2019reexamining}. In contrast to the SPDE route, which is driven by a weak-form discretization of a differential operator, the interpolation route is driven directly by the finite-dimensional distributions of the parent GP. This leads to structured covariance matrices on the grid and is particularly convenient for hierarchical modeling of the bandwidth parameter and the resolution $N$.

We first describe the construction in one dimension. Here we specifically focus on the nodal hat basis $\{\psi_j\}$ in \eqref{eqn: definition of linear basis function}, which satisfies $\psi_j(k/N)=\ind_{\{j=k\}}$. With this basis, the representation
\begin{equation*}
f_N^\star(x) = \sum_{j=0}^N f^\star(j/N)\,\psi_j(x)
\end{equation*}
linearly interpolates a function $f^\star$ between the grid values $\{0,1/N,\ldots,1\}$.
A natural finite-rank GP prior on \eqref{eqn:f_approximation} is then obtained by taking a parent GP prior $f\sim \mbox{GP}(0,K)$ and setting
\begin{equation*}
w^N=\big(f(0), f(1/N), \ldots, f(1)\big)^\top,
\end{equation*}
so that $w^N \sim N(0,G_{N, \theta})$ with $(G_{N, \theta})_{ij}=K(i/N,j/N)$. This gives another instance of the prior in equation \eqref{eqn:w_prior}, with $\theta$ denoting the hyperparameters of the parent GP covariance kernel and $\Sigma_{N,\theta}=G_{N, \theta}$. We refer to inference based on this prior as the {\it GPI (GP interpolation) method}. 

If $K$ is stationary with $K(x,y)=k(x-y)$, then $G_{N, \theta}$ has a Toeplitz structure:
\begin{align}\label{eq: G N}
    G_{N, \theta} =
    \begin{bmatrix}
    k(0) & k\left(\frac{1}{N}\right) & \cdots & k\left(\frac{N}{N}\right)\\
    k\left(\frac{1}{N}\right) & k(0) & \cdots & k\left(\frac{N-1}{N}\right)\\
    \vdots & \vdots & \ddots &  \vdots\\
    k\left(\frac{N}{N}\right) &  k\left(\frac{N-1}{N}\right)  &  \cdots  & k(0)
    \end{bmatrix}.
\end{align}
For $d\ge 2$, one can extend the above construction using a tensor-product hat basis on a regular grid. Define the grid points
\begin{equation*}
    \Big\{\frac{\bj}{N} :=\Big(\frac{j_1}{N}, \ldots, \frac{j_d}{N}\Big): \bj = (j_1, \ldots, j_d)\in\cJ_N\Big\}\subset[0,1]^d.
\end{equation*}
For $f^\star:[0,1]^d\to\bbR$, define
\begin{equation}\label{eqn: definition of f_N d dim}
    f_N^\star(x_1,\ldots,x_d)
    = \sum_{\bj\in\cJ_N}
    w_{\bj}\,\prod_{k=1}^d \psi_{j_k}(x_k),
\end{equation}
with $w_{\bj}=f^\star(\bj/N)$, so that $f_N^\star$ interpolates $f^\star$ at the grid points. As before, one can induce a prior on \eqref{eqn:f_approximation} by taking $w_{\bj}=f(\bj/N)$ for a parent GP $f\sim\mbox{GP}(0,K)$, yielding a multivariate normal prior on $w^N=(w_\bj)_{\bj\in\cJ_N}$. If $K$ is stationary, then $G_{N, \theta}$ again inherits additional structure; for example, for $d=2$, $G_{N, \theta}$ is a block Toeplitz matrix with Toeplitz blocks (BTTB); see, for example, \citet{akaike1973block}. 

\section{Theoretical properties of d-frGPloc}\label{sec:RKHS}

As discussed in Section \ref{sec:2methods}, for fixed $N$, d-frGPloc defines a finite-rank GP prior on $f_N$. 
This section studies two aspects of that prior. First, we describe its reproducing kernel Hilbert space (RKHS), which determines support and the approximation term appearing in posterior contraction arguments. 
Second, we study posterior contraction after placing a prior on the resolution parameter $N$. 
While most of the literature deals with fixed resolution $N$, the common new hierarchy introduced in this paper is the prior on the resolution parameter $N$, which yields adaptive-resolution results for both finite-rank Gaussian process approximations. The treatment of the other hyperparameters follows the corresponding parent-GP theory. 

\subsection{RKHS of a finite-rank GP}

We begin with the RKHS of finite-rank GP prior since it identifies the support and the RKHS norm is directly used in the concentration function defined later. Since the models in Section \ref{sec:2methods} yields Gaussian distribution on the coefficients conditioned on $(N,\theta)$, the RKHS can be written explicitly in terms of the coefficient covariance.

The reproducing kernel Hilbert space (RKHS) determines the support and the geometry of a GP \citep{van2008reproducing}. For a mean-zero GP $f \sim \GP(0, K)$, denote by $\bbH_0$ the set of functions
\begin{align*}
    t \mapsto \sum_{i=1}^k \alpha_i K(s_i,t),
    \qquad s_1,\ldots,s_k\in[0,1]^d,\ \alpha_1,\ldots,\alpha_k\in\bbR,\ k\in\bbN.
\end{align*}
On $\bbH_0$, define an inner product by
\begin{align*}
    \Big\langle \sum_{i=1}^k \alpha_i K(s_i,\cdot), \sum_{j=1}^l \beta_j K(r_j,\cdot) \Big\rangle_{\bbH}
    = \sum_{i=1}^k\sum_{j=1}^l \alpha_i\beta_j K(s_i,r_j),
\end{align*}
and let $\|h\|_{\bbH}^2=\langle h,h\rangle_{\bbH}$ for $h\in\bbH_0$. The RKHS of $f$, denoted by $\bbH$, is defined as the closure of $\bbH_0$ with respect to $\|\cdot\|_{\bbH}$. The closure of $\bbH$ in the ambient Banach space $\bbB$ coincides with the support of the GP \citep{van2008reproducing}. Additional examples of RKHSs associated with different choices of the covariance kernel $K$ can be found in \citet{van2011information}.

Recall that d-frGPloc is defined through the basis expansion \eqref{eqn:f_approximation}, where coefficient vector $w^N=(w_\bj)_{\bj\in\cJ_N}$ follows a Gaussian distribution $w^N \sim N(0, \Sigma_{N, \theta})$ conditioned on $(N,\theta)$. The following lemma describes the RKHS of $f_N$ for fixed $(N,\theta)$.

\begin{lemma}\label{lemma: RKHS of f_N}
Let $f_N=\sum_{\bj\in\cJ_N} w_\bj\psi_\bj$ be the Gaussian process defined on $[0,1]^d$ with linearly independent basis functions $\{\psi_{\bj}: \bj \in \cJ_N\}$ and $w^N\sim N(0,\Sigma_{N, \theta})$, where $\Sigma_{N, \theta}$ is a positive definite matrix indexed by $\cJ_N$. Then the RKHS of $f_N$ is $\mathrm{span}\{\psi_\bj: \bj\in\cJ_N\}$. Moreover, for
$g=\sum_{\bj\in\cJ_N} v_\bj\psi_\bj \in \bbH$ with $v^N=(v_\bj)_{\bj\in\cJ_N}$,
\begin{equation*}
\|g\|_{\bbH}^2=(v^N)^\top \Sigma_{N, \theta}^{-1} v^N.
\end{equation*}
\end{lemma}

Lemma \ref{lemma: RKHS of f_N} shows that for fixed $N$, the support of $f_N$ is the linear span of $\{\psi_j\}$, and the RKHS norm is determined by the coefficient precision matrix $\Sigma_{N,\theta}^{-1}$. In particular, approximation of a globally smooth function depends on both the size of $N$ and the dependence structure in $w^N$.

\subsection{Posterior contraction rates}

We now introduce the posterior contraction framework common to the two constructions. The key new ingredient relative to the parent GP theory is the prior on $N$, and the main focus is how to connect the prior on the coefficient vector $w^N$ with the behavior of the corresponding infinite-dimensional parent GP.

Let $P_2$ denote the prior on $N$ and $\Pi(\cdot\mid N,\theta)$ denote the conditional finite-rank GP prior on $f_N$ given $N$ and $\theta$. When the hyperparameter vector $\theta$ is also assigned a prior (independent of $N$), let $P_1$ denote its probability measure and $p_1$ denote its density. With a prior on $\theta$, the marginal prior on $f_N$ is
\begin{equation*}
    \Pi(\cdot)=\int \Pi(\cdot\mid N,\theta)dP_1(\theta)dP_2(N),
\end{equation*}
and let $\Pi(\cdot\mid\cD_n)$ denote its posterior distribution under \eqref{eqn: fixed regression model}. When $\theta = \theta_n$ is treated as a deterministic sequence depending on $n$ and the smoothness of the true function, $P_1 \equiv \delta_{\theta_n}$ is a Dirac measure at $\theta_n$.

We define $\epsilon_n\to0$ to be a posterior contraction rate of $\Pi(\cdot\mid \cD_n)$ with respect to a metric $d(\cdot,\cdot)$ if there exists $M>0$ such that
\begin{align}\label{eqn: post convergence}
    \bbE_{f^*}\Pi\!\left(f: d(f,f^*)>M\epsilon_n \mid \cD_n\right)\to 0.
\end{align}
A commonly used metric is the empirical $L_2$ norm $\|h\|_n^2=n^{-1}\sum_{i=1}^n h(x_i)^2$ \citep{Van2008Rates}. For regular fixed designs, the minimax rate over $C^\alpha[0,1]^d$ in this metric is $n^{-\alpha/(2\alpha+d)}$ \citep{tsybakov2008introduction}.

The concentration function is a useful tool to study posterior concentration in standard statistical models with Gaussian process priors \citep{van2008reproducing}. Since our model is Gaussian conditional on $(N,\theta)$, we consider the conditional concentration function
\begin{align}\label{eqn: def of the concentration function}
    \phi_{f^*}^{(N,\theta)}(\epsilon)
    =
    \inf_{h\in\bbH_{N,\theta}:\ \|h-f^*\|_\infty\le \epsilon}
    \frac12\|h\|_{\bbH_{N,\theta}}^2
    \;-\;
    \log \Pi\!\left(\|f_N\|_\infty<\epsilon \mid N,\theta\right),
\end{align}
where $\bbH_{N,\theta}$ is the RKHS of $f_N$ under $\Pi(\cdot\mid N,\theta)$. This quantity controls the conditional prior mass around $f^*$ up to constants \citep{van2008reproducing}. Under the marginal prior,
\begin{equation*}
    \Pi\big(\|f-f^*\|_\infty<\epsilon\big)
=
\int \Pi\big(\|f_N-f^*\|_\infty<\epsilon \mid N,\theta\big)\,p_1(\theta)\,d\theta\,dP_2(N).
\end{equation*}

The required resolution $N$ to make the prior mass $\Pi\big(\|f_N-f^*\|_\infty<\epsilon \mid N,\theta\big)$ positive depends on the H\"older smoothness of $f^*$.  For $\alpha>0$, let
$\underline{\alpha}$ denote the largest integer strictly smaller than $\alpha$, and
for a multi-index $\bk=(k_1,\dots,k_d)\in\bbN_0^d$ write $|\bk|=k_1+\cdots+k_d$ for
its order and $D^\bk=\partial^{|\bk|}/(\partial x_1^{k_1}\cdots\partial x_d^{k_d})$ for the
corresponding mixed partial derivative. The H\"older space $C^\alpha[0,1]^d$ is the
set of functions $f:[0,1]^d\to\bbR$ whose derivatives $D^\bk f$ exist for all
$|\bk|\le\underline{\alpha}$ and for which
\begin{equation}\label{eqn: holder norm}
    \|f\|_{C^\alpha}
    := \max_{\bk:\,|\bk|\le\underline{\alpha}}\,\sup_{x\in[0,1]^d}|D^\bk f(x)|
     + \max_{\bk:\,|\bk|=\underline{\alpha}}\,\sup_{x\neq y}
       \frac{|D^\bk f(x)-D^\bk f(y)|}{\|x-y\|_2^{\,\alpha-\underline{\alpha}}}
     <\infty.
\end{equation}
The following lemma gives a resolution level, as a function of $\alpha$, that
suffices to approximate $f^*$ within supremum-norm error $\epsilon$ using the
tensor-product hat basis defined in Section~\ref{ssec:GPI}.
\begin{lemma}\label{lemma: necessary condition for N}
Assume $\{\psi_j\}_{j\in\cJ_N}$ is the tensor-product hat basis on the regular grid $\{j/N: j\in\cJ_N\}\subset[0,1]^d$. For all $f^* \in C^{\alpha}[0,1]^d$, there exists $N\asymp \epsilon^{-1/(\alpha\land 2)}$ such that
\begin{equation*}
\Big\{ f_N \in \mathrm{span}\{\psi_j: j\in\cJ_N\}:\ \|f_N-f^*\|_\infty<\epsilon \Big\}\neq\emptyset.
\end{equation*}
\end{lemma}

Lemma \ref{lemma: necessary condition for N} implies that for the approximation at accuracy $\epsilon$, the resolution $N$ to be sufficiently large, with the order depending on the smoothness of $f^*$. This motivates assigning a prior on $N$ to achieve adaptivity. 

\subsection{Posterior contraction rate of the SPDE approach}\label{ssec:SPDEt}

We first consider the posterior contraction rate result of the SPDE method. The main theorem in this subsection implies that with a prior placed on $N$, the finite-rank SPDE prior inherits the same contraction rate as the parent Mat\'ern GP, under the same sample-size- and smoothness-dependent bandwidth scaling. 

The main technical difficulty is to connect the parent GP under a Neumann boundary condition with the finite rank SPDE prior. The property of the parent GP is analyzed based on spectral decomposition and is compared with the Gaussian prior on the coefficient vectors as suggested in \citet{sanz2022finite}, while the parent-GP contraction result is taken from \citet{fang2025posterior}.

As mentioned in Section \ref{sec:2methods}, $f_N$ derived from the SPDE method is defined on the bounded domain $D=[0,1]^d$ with homogeneous Neumann boundary conditions. The operator $\cL_\kappa=\kappa^2\mbox{id}-\Delta$ admits an orthonormal eigenbasis indexed by $\bk=(k_1,\ldots,k_d)\in\bbN_0^d$,
\begin{equation*}
    e_{\bk}(x)=\prod_{m=1}^d e_{k_m}(x_m),
    \qquad
    e_{0}(t)=1,\quad e_{\ell}(t)=\sqrt{2}\cos(\ell\pi t)\ (\ell\ge1),
\end{equation*}
and
\begin{equation*}
\cL_\kappa e_{\bk}=\{\kappa^2+\pi^2\|\bk\|_2^2\}e_{\bk}.
\end{equation*}
We use this representation to define a Besov-type smoothness class adapted to Neumann boundary conditions, following Section~3 of \cite{sanz2022finite}. This class measures smoothness through how quickly the low-frequency components approximate the function in supremum norm. Let $S$ be an even function on $\bbR$ satisfying
\begin{equation*}
0\le S\le 1,\quad S(x)=1\ \text{for }x\in[-0.5,0.5],\quad S(x)=0\ \text{for }|x|>1,
\end{equation*}
and define $S_j(x)=S(2^{-j}x)$. For $f=\sum_{\bk\in\bbN_0^d} f_{\bk} e_{\bk}$, define
\begin{equation*}
    \TS_j f = \sum_{\bk\in\bbN_0^d} S_j(\pi\|\bk\|_2)\, f_{\bk} e_{\bk}
\end{equation*}
and
\begin{equation*}
    \Besova = \Big\{ f=\sum_{\bk\in\bbN_0^d} f_{\bk} e_{\bk}:\ 
\|f\|_{\Besova} := \sup_{j\in\bbN} 2^{\alpha j}\|\TS_j f - f\|_\infty <\infty \Big\}.
\end{equation*}
This is a Neumann-adapted Besov-type function class with smoothness $\alpha$. 

Theorem 3.4 of \citet{fang2025posterior} shows that for the parent Mat\'ern GP with covariance kernel \eqref{eqn: matern cov kernel}, if the true function has H\"older and Sobolev smoothness $\alpha$ and
\begin{equation}\label{eq:oracle kappa in SPDE}
    \kappa_n = n^{\frac{\beta - d/2 - \alpha}{(2\alpha + d)(\beta - d/2)}},
\end{equation}
then under $\alpha < \beta - d/2$ the posterior contracts at the minimax rate $n^{-\alpha/(2\alpha+d)}$. The theorem below shows that, under the same bandwidth scaling, adding a prior on the resolution $N$ preserves this rate for the finite-rank SPDE prior. Thus the SPDE result is adaptive in $N$, with the bandwidth parameter inherited from the parent Mat\'ern theory.

\begin{theorem}
(Posterior contraction for the SPDE method)
\label{thm:SPDE optimality}
Let $0<\alpha<\beta-d/2$, let $f^* \in \Besova$, and let $f_N$ be a d-frGPloc defined in \S \ref{ssec:SPDE}. Suppose $\kappa_n$ satisfies \eqref{eq:oracle kappa in SPDE} and the prior on $N$ satisfies
\begin{equation}
    P_2(N \ge x) \ge K_1 \exp\!\big(-K_2 (\log x)^{K_3}\big), \quad \forall x>x_0, \label{eqn: prior on N SPDE}
\end{equation}
for some constants $x_0,K_1,K_2,K_3>0$. Then \eqref{eqn: post convergence} holds with $\epsilon_n=n^{-\alpha/(2\alpha+d)}(\log n)^C$ for some $C>0$.
\end{theorem}
The tail condition on $P_2$ is mild; for example, any polynomial prior $P_2(N=m)\propto m^{-a}$ with $a>1$, or log-normal prior on $N$ satisfies \eqref{eqn: prior on N SPDE}.

The proof is based on a spectral comparison between the continuous Neumann operator $\cL_\kappa$ and its mass-lumped finite-element discretization. Evaluating the Neumann eigenfunctions $e_{\bk}$ at the grid points yields vectors that are generalized eigenvectors of the discrete precision matrix $Q$ with respect to the lumped mass matrix $C$; the resulting discrete eigenvalues are comparable to the continuous ones up to a factor depending only on $\beta$ (Lemma \ref{lem:SPDE_eigs_d_app} in the Supplementary Material). This equivalence implies that the RKHS norm and centered small-ball probability of the finite-rank prior match those of the parent Mat\'ern GP up to constants for large enough $N$. From this, the three conditions of the posterior contraction theorem of \citet{Van2008Rates} (prior concentration, sieve probability, and metric entropy) are verified at the same rate (up to logarithmic factors) as for the parent GP. The lower tail condition \eqref{eqn: prior on N SPDE} ensures that such resolutions receive enough prior mass. Additionally, a single sieve can be utilized to cover all sufficiently large $N$, leading to the same order of metric entropy up to logarithmic factor. The finite-element setup follows the bounded-domain Neumann framework of \citet{sanz2022finite}, while the parent Mat\'ern benchmark is adapted from \citet{fang2025posterior}. 

\subsection{Posterior contraction rate of the GP interpolation approach}\label{ssec:GPIt}

We next turn to the interpolation construction. 
Unlike the SPDE result above, the GPI result assigns priors to both the bandwidth parameter $\kappa$ and the resolution parameter $N$. Thus the interpolation construction retains the adaptive bandwidth behavior of the parent squared-exponential GP while also learning the finite-rank resolution.
The argument combines the parent squared-exponential GP concentration result with an interpolation error bound showing that, once $N$ is large enough relative to $\kappa$, the finite-rank approximation is accurate at the relevant scale.

Motivated by \citet{van2009adaptive}, we consider a squared-exponential covariance kernel with inverse-bandwidth parameter $\kappa$,
\begin{align}\label{eqn: RBF kernel function}
    K(x_1,x_2)=\exp\big(-\kappa^2\|x_1-x_2\|_2^2\big).
\end{align}
\citet{van2009adaptive} showed that when the original GP is used as a prior for $f^*$, the minimax optimal posterior contraction rate can be achieved by placing a gamma-type prior on $\kappa$. 

We show that adding a prior on $N$ to the interpolation construction preserves this parent-GP adaptive bandwidth behavior, while allowing the adaptive finite-rank resolutions.

\begin{theorem}(Adaptive posterior contraction of GPI)\label{thm:GPI optimality}
Let $f^* \in C^{\alpha}[0, 1]^d$ with $\alpha > 0$ and $f_N$ be a d-frGPloc defined in \S \ref{ssec:GPI} with parent GP covariance \eqref{eqn: RBF kernel function}, and assume independent priors on $\kappa$ and $N$ satisfying
\begin{align}
    C_1 \kappa^p \exp\!\big(-D_1 \kappa^d (\log \kappa)^q\big)
    &\le p_1(\kappa) \le
    C_2 \kappa^p \exp\!\big(-D_2 \kappa^d (\log \kappa)^q\big), \label{eqn: prior on kappa GPI}\\
    P_2(N \ge x) &\ge K_1 \exp\!\big(-K_2 (\log x)^{K_3}\big), \quad \forall x>x_0 \label{eqn: prior on N GPI}
\end{align}
for some $C_1, C_2, D_1, D_2, p, q, K_1, K_2, K_3, x_0 > 0$.
Then $\epsilon_n=n^{-\alpha/(2\alpha+d)}(\log n)^C$ for some $C>0$ satisfies \eqref{eqn: post convergence}.
\end{theorem}
The condition on $p_1$ is the gamma-type inverse-bandwidth prior used in \citet{van2009adaptive}; the constants control the polynomial prefactor and logarithmic correction and affect only the logarithmic factor in the contraction rate.

The proof of Theorem \ref{thm:GPI optimality} makes use of the optimal prior concentration of a GP with the squared-exponential kernel \citep{van2009adaptive}. We select $N$ large enough, depending on $\kappa$, so that the interpolation error $f_N-f$ is negligible at the scale relevant for prior concentration, where $f$ denotes the parent GP. Combining this with \eqref{eqn: prior on N GPI} yields the desired prior concentration, and a single sieve can be utilized to cover all sufficiently large $N$, leading to the same order of metric entropy up to logarithmic factor. See Appendix \ref{appendix:GPI} of the Supplementary Material. 

\section{Posterior computation for the adaptive resolution}\label{sec:comp}

In practice, the adaptive-resolution formulation is computationally inconvenient because the dimension of the coefficient vector $w^N$ changes with $N$.
We introduce an efficient computational scheme below, where we avoid reversible-jump updates by first sampling $N$ and the kernel hyperparameters from the {\it marginal posterior} obtained after integrating out $w^N$, and then sampling $w^N$ from its Gaussian conditional posterior.
This method gives an additional benefit for the SPDE construction, because the marginal posterior of $N$ and $\theta$ can be evaluated efficiently using the sparse precision matrix of $w^N$ conditional on $N$ and $\theta$.

Recall from Section \ref{ssec:GPI} the finite-rank representation
\begin{equation}\label{eq:comp_fN}
    f_N(x)=\sum_{\bj\in\cJ_N} w_\bj\psi_\bj(x),
    \qquad
    \cJ_N=\{\bj=(j_1,\ldots,j_d): j_k\in\{0,1,\ldots,N\}\},
\end{equation}
and let $m_N=|\cJ_N|=(N+1)^d$.
For both GPI and SPDE methods, a prior distribution on the coefficient vector $w^N=(w_j)_{j\in\cJ_N}$ given the resolution $N$ and the hyperparameter $\theta$ can be written as
\begin{equation*}
    w^N \mid N, \theta \sim N(0, \Sigma_{N, \theta}),
\end{equation*}

For GPI, $\theta=\kappa$ and $\Sigma_{N,\theta}=G_{N, \theta}$, where $G_{N, \kappa}$ is the parent-GP covariance matrix on the grid. For SPDE, $\theta=\kappa$ and $\Sigma_{N,\theta}=Q_{N,\theta}^{-1}$, where $Q_{N,\theta}$ is the mass-lumped finite-element precision matrix.
Since $\theta = \kappa$ for both constructions, we write $\kappa$ in place of $\theta$ for the remainder of this section.

We place independent priors $N\sim p(N)$ and $\kappa\sim p_1(\kappa)$. For the theoretical SPDE result in Theorem \ref{thm:SPDE optimality}, $\kappa$ is fixed at the sample-size-dependent value $\kappa_n$. In that case, the update below is applied only to $N$. When a hierarchical implementation is used, as in the simulations, the marginal posterior of $(N,\kappa)$ is calculated.

Let $y=(y_1,\ldots,y_n)^\top$ and define the design matrix $\Phi_N\in\bbR^{n\times m_N}$ by
\begin{equation}\label{eq:comp_Phi}
    (\Phi_N)_{i,j}=\psi_j(x_i),\qquad i=1,\ldots,n,\ j\in\cJ_N.
\end{equation}
Then $(f_N(x_1),\ldots,f_N(x_n))^\top=\Phi_N w^N$ holds and the likelihood under \eqref{eqn: fixed regression model} is written as
\begin{equation}\label{eq:comp_like}
    p(y\mid w^N,N,\kappa)\propto
    \exp\Big\{-\frac{1}{2\sigma^2}(y-\Phi_N w^N)^\top(y-\Phi_N w^N)\Big\}.
\end{equation}
Since $\{\psi_j\}$ are locally supported, each row of $\Phi_N$ has at most $2^d$ nonzero entries, so $\Phi_N^\top y$ and $\Phi_N^\top\Phi_N$ can be formed in $O(n)$ operations. By conjugacy, the conditional posterior of $w^N$ given $(N,\kappa)$ is Gaussian:
\begin{equation}\label{eq:comp_post_w}
    w^N \mid N,\kappa,\Data \sim N(\mu_{N,\kappa},\Sigma^*_{N,\kappa}),
\end{equation}
where
\begin{equation}\label{eq:comp_post_params}
    (\Sigma^*_{N,\kappa})^{-1}=\Sigma_{N,\kappa}^{-1}+\sigma^{-2}\Phi_N^\top\Phi_N,
    \qquad
    \mu_{N,\kappa}=\sigma^{-2}\Sigma^*_{N,\kappa}\Phi_N^\top y.
\end{equation}

Therefore, it remains to sample $(N, \kappa)$ from the marginal posterior. From
\begin{equation}\label{eq:comp_post_Nkappa}
    p(N,\kappa\mid\Data)\propto p(N)p_1(\kappa)p(\Data\mid N,\kappa) = p(N)p_1(\kappa) \int p(w^N, \Data \mid N, \kappa) dw^N,
\end{equation}
we obtain
\begin{equation}\label{eq:comp_marg_like}
    \log p(\Data\mid N,\kappa)
    =
    C_y
    -\frac12\log|\Sigma_{N,\kappa}|
    +\frac12\log|\Sigma^*_{N,\kappa}|
    +\frac{1}{2\sigma^4} y^\top \Phi_N  \Sigma^*_{N,\kappa} \Phi_N^T y,
\end{equation}
where $C_y$ does not depend on $(N,\kappa)$.
From this, we update $N$ using a local proposal on its discrete support and update $\kappa$ using a Metropolis--Hastings proposal density $q_{\kappa}(\cdot\mid\kappa)$ on $(0,\infty)$.
Given a current state $(N,\kappa)$, propose $(N',\kappa')$ and accept with probability
\begin{equation}\label{eq:comp_accept}
    \min\Big\{1,\
    \frac{p(N')p_1(\kappa')p(\Data\mid N',\kappa')}{p(N)p_1(\kappa)p(\Data\mid N,\kappa)}
    \cdot
    \frac{q_N(N\mid N')q_{\kappa}(\kappa\mid\kappa')}{q_N(N'\mid N)q_{\kappa}(\kappa'\mid\kappa)}
    \Big\},
\end{equation}
where $p(\Data\mid N,\kappa)$ is evaluated via \eqref{eq:comp_marg_like}.
After updating $(N,\kappa)$, we draw $w^N$ from \eqref{eq:comp_post_w}.

The dominant computation per iteration is governed by the coefficient dimension $m_N$ rather than the sample size $n$:
forming $\Phi_N^T y$ and $\Phi_N^\top \Phi_N$ costs $O(n)$, while dense linear algebra for evaluating \eqref{eq:comp_marg_like} and sampling \eqref{eq:comp_post_w} costs $O(m_N^3)$.
Thus the procedure is computationally attractive when $m_N\ll n$, even though dense factorizations remain the bottleneck as $d$ and $N$ grow.
When the parent kernel is stationary and the grid is regular, $\Sigma_{N,\kappa}$ has Toeplitz structure for $d=1$ and BTTB structure for $d=2$, which can be exploited to accelerate linear algebra operations involving $\Sigma_{N,\kappa}$, although not for $\Sigma_{N, \kappa}^*$ \citep{ray2020efficient}.
For the SPDE method, however, both $\Sigma^{-1}_{N,\kappa}$ and $(\Sigma^*_{N,\kappa})^{-1}$ remain sparse. This is because of the local support of the basis functions, leading to the sparsity of $\Sigma^{-1}_{N,\kappa}$ and $\Phi_N^\top\Phi_N$.

Therefore, utilizing sparse Cholesky factorization, the SPDE method avoids dense $m_N \times m_N$ factorizations. In one dimension, where the precision matrix is banded, the coefficient-space linear algebra is linear in $m_N$ up to constants depending on the bandwidth. In higher dimensions, the exact cost depends on the sparsity in Cholesky factorization, which requires at most $O(m_N^3)$. Thus a direct implementation of GPI has cost $O(n+m_N^3)$ per iteration. For the SPDE method, the Cholesky factorization of $\big(\Sigma_{N, \kappa}^*\big)^{-1}$ can be obtained more efficiently thanks to its sparsity. For instance, in one-dimensional case, the Cholesky decomposition of $\big(\Sigma_{N, \kappa}^*\big)^{-1}$ requires $O(m_N)$ as it has bandwidth $\beta$, leading to per iteration computational cost $O(n + m_N)$.
Therefore, both methods remain scalable for larger $n$, by avoiding dense $n \times n$ parent-GP covariance matrix operations.

{\normalsize
\begin{algorithm}[htbp!]
\caption{Hierarchical sampling with continuous $\kappa$}
\label{alg:GPI_MH_continuous}
\KwIn{Observed data $\Data=\{(x_i,y_i)\}_{i=1}^n$, priors $p(N)$ and $p_1(\kappa)$}
\KwOut{Posterior samples $\{(w^{N},N,\kappa)^{(t)}\}_{t=1}^T$}
Initialize $(N^{(0)},\kappa^{(0)})$\;

\For{$t=1,\ldots,T$}{
  Propose $N'\sim q_N(\cdot\mid N^{(t-1)})$ and $\kappa'\sim q_{\kappa}(\cdot\mid \kappa^{(t-1)})$\;

  Evaluate $\log p(\Data\mid N',\kappa')$ and $\log p(\Data\mid N^{(t-1)},\kappa^{(t-1)})$ via \eqref{eq:comp_marg_like}\;

  Accept $(N',\kappa')$ with probability \eqref{eq:comp_accept}; otherwise keep $(N^{(t)},\kappa^{(t)})=(N^{(t-1)},\kappa^{(t-1)})$\;

  Sample $w^{N^{(t)}} \sim N(\mu_{N^{(t)},\kappa^{(t)}},\Sigma^*_{N^{(t)},\kappa^{(t)}})$ using \eqref{eq:comp_post_w}--\eqref{eq:comp_post_params}\;
}
\end{algorithm}
}

\section{Simulation studies}\label{sec:sims}

We evaluate the empirical performance of the two finite-rank approximations introduced in Section \ref{sec:2methods} under the fixed-design regression model \eqref{eqn: fixed regression model}.
The simulations are designed to examine two questions.
First, how closely do the adaptive finite-rank procedures track the corresponding parent GP procedures?
Second, how much is gained by learning the resolution parameter $N$ rather than fixing it in advance?

For each experiment we generate design points $x_i$ uniformly on $[0,1]^d$ and sample
\begin{equation}\label{eq:sims_model}
    y_i=f(x_i)+\epsilon_i,\qquad \epsilon_i\sim N(0,\sigma^2),
\end{equation}
with $\sigma^2=0.01$ treated as known.
As a point estimator we use the posterior mean
\begin{equation}\label{eq:posterior_mean_sim}
    \hat f(x)=\bbE\{f(x)\mid\Data\},
\end{equation}
approximated by Monte Carlo averages of the retained posterior draws.
Performance is measured by the average mean squared error (AMSE) on a regular evaluation grid:
\begin{equation}\label{eq:amse_def_rewrite}
    \mathrm{AMSE}(\hat f)
    =
    \frac{1}{|\cG_K|}\sum_{t\in\cG_K}\{\hat f(t)-f(t)\}^2,
\end{equation}
where $\cG_K$ denotes the $K^d$ regular grid on $[0,1]^d$.
In one dimension we take $K=1000$, and in two dimensions we take a $25\times 25$ grid.

For the SPDE method, Theorem \ref{thm:SPDE optimality} covers adaptive resolution under the Mat\'{e}rn bandwidth scaling depending on the true function smoothness. In the simulations, we report an empirical hierarchical version that updates $\kappa$ jointly with $N$. This implementation is used to examine practical sensitivity to bandwidth learning, while Theorem \ref{thm:SPDE optimality} deals with the fixed $\kappa$ scaling depending on the true function smoothness.

\subsection{Truth functions and simulation setup}

In $d=1$ we consider the Fourier-series family
\begin{equation}\label{eq:fourier_truth_rewrite}
    f_{\alpha}(x)
    =
    \sum_{j=1}^{500}\sqrt{2}\,\sin(j)\,\cos\big(\pi(j-\tfrac12)x\big)\,j^{-(1+\alpha)},
    \qquad x\in[0,1],
\end{equation}
with $\alpha\in\{0.7,\,2.5\}$.
These represent a relatively rough truth and a smoother truth.
In $d=2$ we consider
\begin{equation}\label{eq:truth_2d}
    f_3(x_1,x_2)=\sin(5|x_1-0.7|+2x_2)+2x_2^2,\qquad (x_1,x_2)\in[0,1]^2.
\end{equation}

For both GPI and SPDE we use piecewise linear basis functions in $d=1$ and tensor-product piecewise linear bases in $d=2$. In the adaptive finite-rank procedures, we place a prior on $N$ and a continuous prior on $\kappa$. Specifically, in the simulation study we take
\begin{equation}\label{eq:priors_cont_sim}
p(N)\propto N^{-2},
\qquad
p_1(\kappa)\propto \kappa^{-1}\ind\{\underline{\kappa}\le \kappa \le \overline{\kappa}\},
\end{equation}
on prespecified supports. The prior $p(N)\propto N^{-2}$ satisfies the lower-tail condition in Theorems \ref{thm:SPDE optimality} and \ref{thm:GPI optimality}, since $\sum_{m\ge x}m^{-2}\asymp x^{-1}$. Thus $N$ is discrete, while $\kappa$ is continuous on a bounded interval.

For GPI, posterior sampling is carried out by Algorithm \ref{alg:GPI_MH_continuous}. For SPDE, we use the analogous hierarchical update over $(N,\kappa)$ together with Gaussian coefficient sampling conditional on $(N,\kappa)$. For each adaptive finite-rank procedure, we ran 5000 MCMC iterations, discarded the first 2500 as burn-in, and retained the remaining 2500 posterior draws for posterior means and posterior summaries. Each simulation configuration was repeated over 20 independent replications.

In one dimension we compare each adaptive finite-rank procedure with its corresponding parent GP procedure.
For GPI, the parent prior is the squared-exponential GP with covariance kernel \eqref{eqn: RBF kernel function}.
For SPDE, the parent prior is the Mat\'ern GP corresponding to \eqref{eqn: matern cov kernel} with the same smoothness specification.
In two dimensions we focus on the finite-rank procedures and compare the adaptive method with fixed-$N$ baselines, since exact parent-GP computation is substantially more expensive in the sample size range considered here.

\subsection{Comparison with the parent GP in one dimension}\label{ssec:sims_rate}

\begin{figure}[tb]
  \centering
  \includegraphics[width=0.9\linewidth]{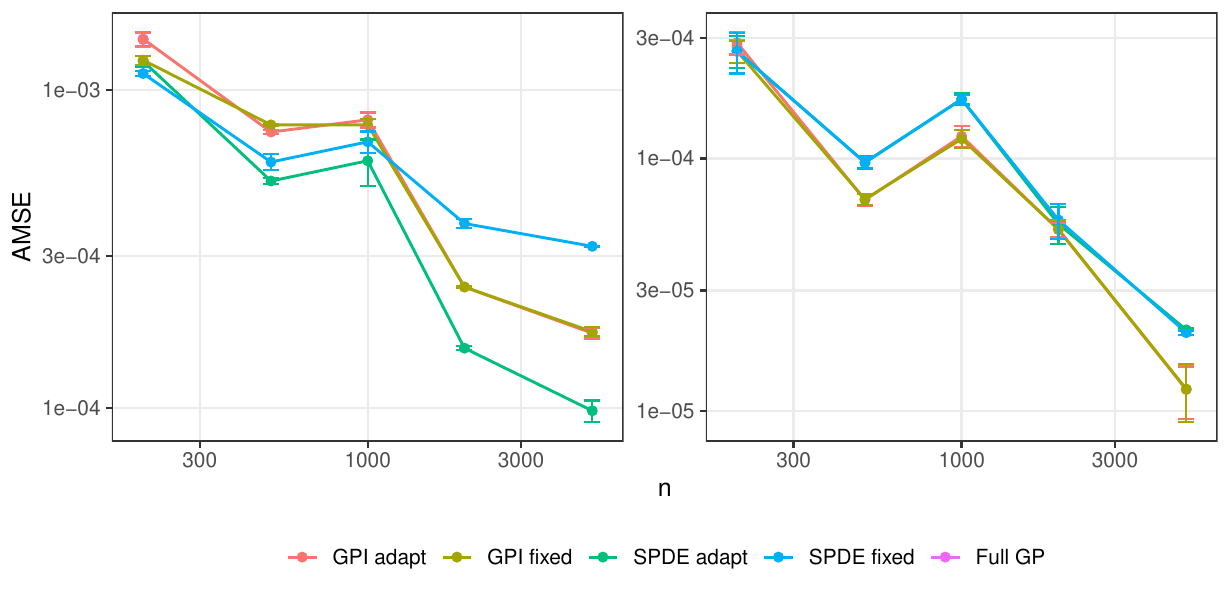}
  \caption{One-dimensional AMSE versus $n$ for the Fourier-series truths \eqref{eq:fourier_truth_rewrite} with $\alpha=0.7$ (left) and $\alpha=2.5$ (right). Curves compare adaptive and fixed-resolution finite-rank procedures with their corresponding parent GP baselines. Points denote the median over 20 independent replications and error bars denote the interquartile range.}
  \label{fig:rate_1d}
\end{figure}

Figure \ref{fig:rate_1d} reports AMSE as a function of $n$ in $d=1$ for $\alpha\in\{0.7,2.5\}$.
Points show the median AMSE over 20 independent replications and error bars show the interquartile range.
For each truth, the adaptive GPI procedure is compared with the exact squared-exponential GP baseline, and the adaptive SPDE procedure is compared with the Mat\'ern GP baseline.

In both smoothness regimes, the adaptive finite-rank procedures closely track the predictive performance of their corresponding parent GP procedures.
For GPI, this agrees with the hierarchical theory over both $N$ and $\kappa$. For SPDE, the simulation uses a hierarchical implementation over $(N,\kappa)$, while the theorem establishes adaptive resolution under the parent Mat\'ern bandwidth scaling.

\subsection{The effect of adaptive resolution in two dimensions}\label{ssec:sims_rate_2d}

\begin{figure}[tb]
  \centering
  \includegraphics[width=0.9\linewidth]{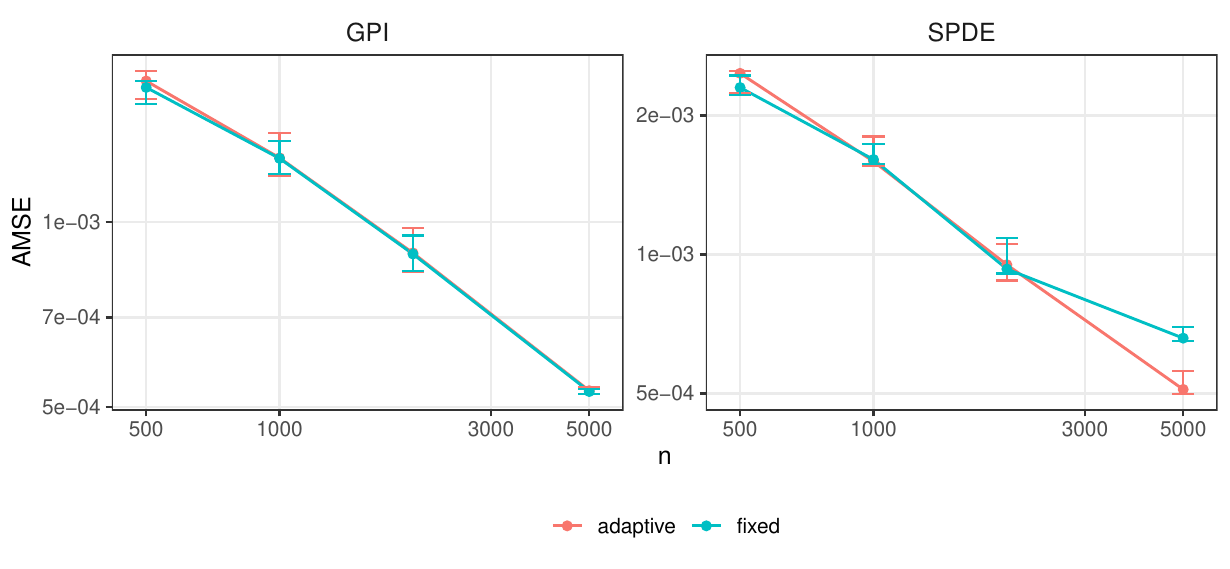}
  \caption{Two-dimensional AMSE versus $n$ for $f_3$. The adaptive-$N$ procedure is compared with fixed-resolution baselines for GPI and SPDE. Points denote the median over 20 independent replications and error bars denote the interquartile range.}
  \label{fig:rate_2d}
\end{figure}

To isolate the role of the resolution parameter, Figure \ref{fig:rate_2d} compares the adaptive procedure with several fixed-$N$ baselines in $d=2$ for the truth \eqref{eq:truth_2d}.
The fixed-$N$ curves show clear sensitivity to the chosen grid resolution.
When $N$ is too small, the approximation is overly restrictive; when $N$ is too large, computation becomes more expensive without a uniform gain in finite-sample performance.
In contrast, the adaptive procedure tracks the better fixed-$N$ choices across sample sizes without manual tuning.

This behavior is consistent with the theory.
For fixed $N$, the prior is supported on a restricted linear span, so approximation quality depends directly on the chosen resolution.
By placing a prior on $N$, the posterior can adapt the effective model complexity to the data.

\subsection{Posterior adaptation of resolution and hyperparameters}\label{ssec:sims_adapt_diag}

\begin{figure}[tb]
  \centering
  \includegraphics[width=0.95\linewidth]{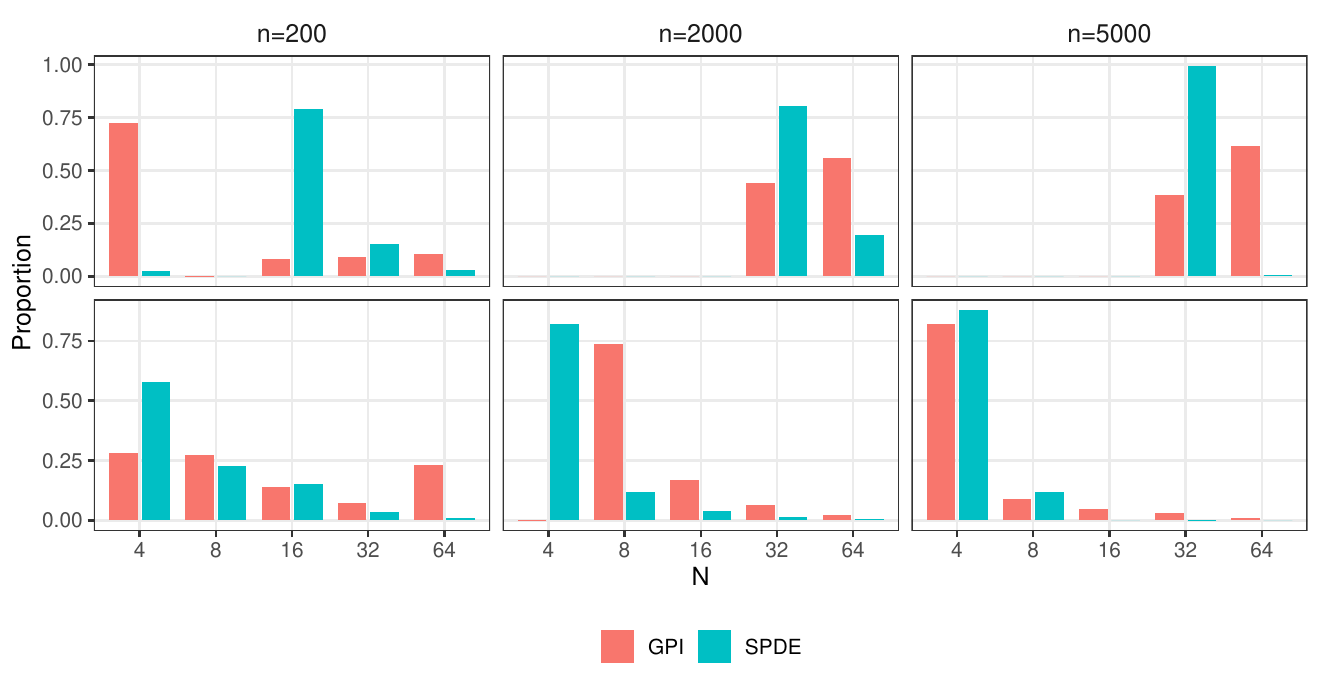}
  \caption{Posterior samples of $N$ in $d=1$ under the adaptive hierarchy, shown for different $n$, under true functions with $\alpha=0.7$ (top) and $\alpha=2.5$ (bottom).}
  \label{fig:postN_1d}
\end{figure}

Figure \ref{fig:postN_1d} summarizes posterior draws of $N$ in $d=1$ across $n$ and across the two smoothness regimes.
For the rougher truth, the posterior places more mass on larger resolutions as $n$ increases, reflecting the need to represent higher-frequency variation.
For the smoother truth, the posterior concentrates on smaller values of $N$, indicating that finer resolutions offer less additional benefit.

\begin{figure}[tb]
  \centering
  \includegraphics[width=0.95\linewidth]{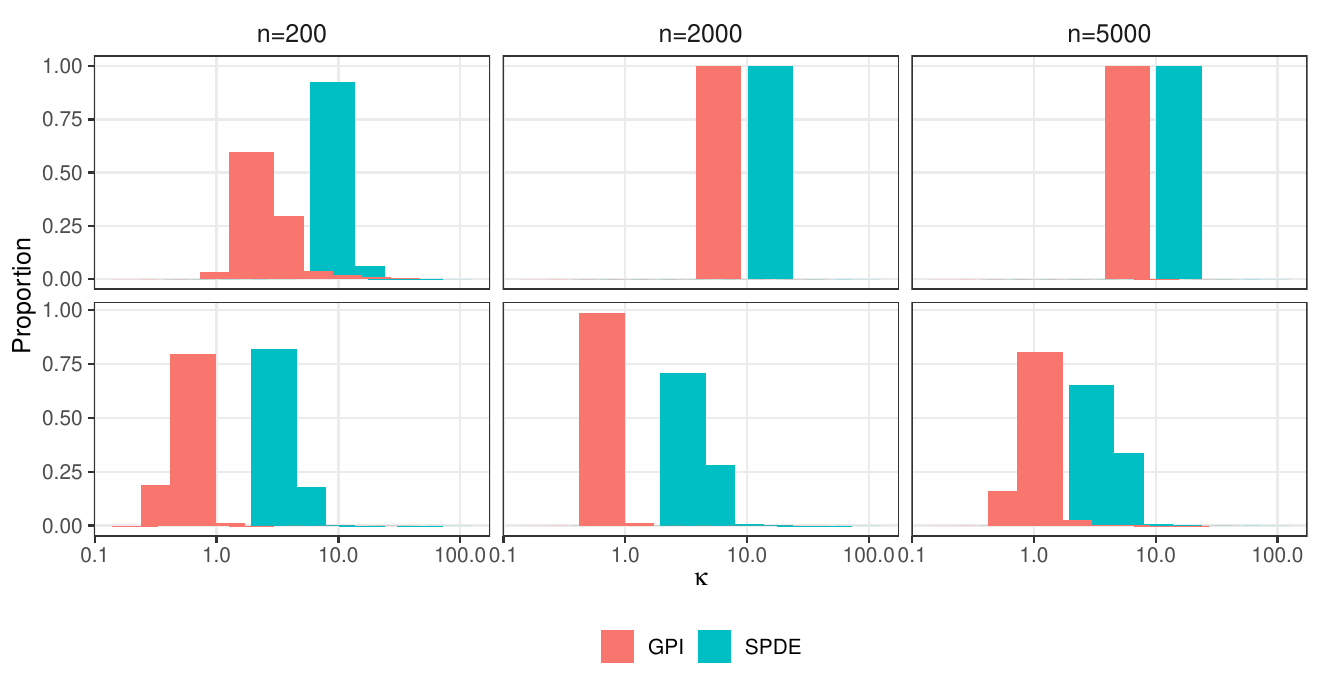}
  \caption{Posterior samples of $\kappa$ in $d=1$ under the adaptive hierarchy, shown for different $n$, under true functions with $\alpha=0.7$ (top) and $\alpha=2.5$ (bottom).}
  \label{fig:postKappa_1d}
\end{figure}

Figure \ref{fig:postKappa_1d} reports posterior draws of $\kappa$.
Since $\kappa$ is treated as a continuous parameter, the posterior is summarized through its sampled distribution rather than through exact-value posterior masses.
The posterior distribution of $\kappa$ shifts with $n$ and differs across the two smoothness regimes, showing that the scale parameter adapts jointly with $N$ under the hierarchical specification.

\subsection{Runtime scaling}\label{ssec:sims_runtime}

\begin{figure}[tb]
  \centering
  \includegraphics[width=0.9\linewidth]{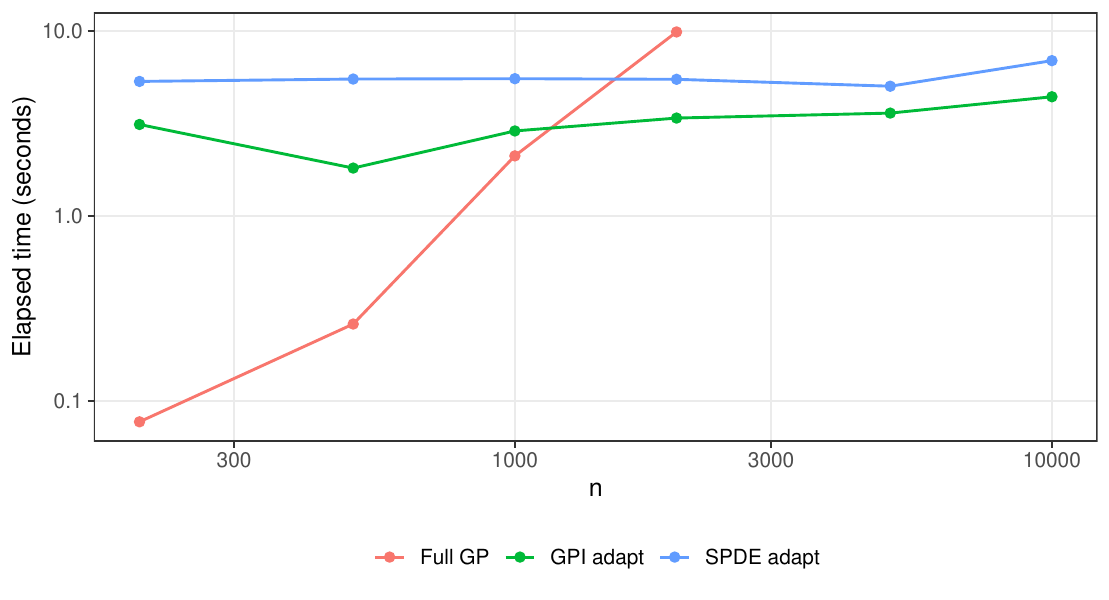}
  \caption{Runtime (seconds) versus $n$ in $d=1$. The parent GP baselines become rapidly more expensive with $n$, while GPI and SPDE remain comparatively stable over the same range. Points denote the median over 20 independent replications and error bars denote the interquartile range.}
  \label{fig:runtime_1d}
\end{figure}

Finally, Figure \ref{fig:runtime_1d} reports wall-clock runtime in $d=1$ as a function of $n$.
The parent GP baselines require factorizations of $n\times n$ covariance matrices and become rapidly more expensive as $n$ grows.
In contrast, the dominant linear algebra for the finite-rank methods is governed by the coefficient dimension $m_N=(N+1)^d$ rather than $n$.
Over the range considered here, the adaptive GPI and adaptive SPDE procedures remain substantially cheaper while maintaining predictive accuracy comparable to that of the corresponding parent GP procedures.

\section{Discussion}\label{sec:disc}

This paper establishes that finite-rank Gaussian process approximations with locally supported bases and dependent Gaussian coefficients can inherit the posterior contraction rates of their corresponding parent GP priors when a suitable prior is placed on the resolution parameter $N$. We have shown this for two representative constructions: the finite-element SPDE approximation of Mat\'ern GPs and the lattice-based interpolation of squared-exponential GPs. The common methodological contribution is \textit{adaptive resolution}: the finite-rank dimension $|{\mathcal J}_N| = (N+1)^d$ is learned from data rather than fixed in advance, preserving the parent-GP contraction rate while maintaining computational advantages from the locally supported basis structure.

For the SPDE approach (Theorem~3.3), adaptation is over the resolution $N$ under the bandwidth scaling $\kappa_n$ inherited from the parent Mat\'ern theory \citep{fang2025posterior}. The proof relies on spectral comparison between the continuous Neumann eigenbasis and the mass-lumped finite-element discretization, showing that the finite-rank RKHS norm and small-ball probability match those of the parent GP up to constants for sufficiently large $N$. This connects the numerical-approximation perspective of \citet{sanz2022finite} with the statistical contraction theory of Mat\'ern GPs.

For the GPI approach (Theorem~3.4), adaptation is over both the bandwidth $\kappa$ and the resolution $N$. The hierarchical specification with horseshoe or gamma-type priors on $\kappa$ recovers the minimax-optimal rate up to logarithmic factors, matching the fully adaptive parent squared-exponential GP result of \citet{van2009adaptive}. The proof uses an interpolation error bound to show that $\|f_N - f\|_\infty$ is negligible when $N$ is large enough relative to $\kappa$, allowing the finite-rank prior concentration to inherit the parent-GP behavior. This strategy has a conceptual parallel with the oracle-GP comparison approach of \citet{castillo2025deep}, where an oracle GP serves as an intermediate reference for establishing contraction rates.

\textbf{Extensions and open problems.} Several natural directions remain open:


\textit{Higher dimensions.} The $O((N+1)^{3d})$ cost makes direct application challenging for $d \geq 5$. Structured approaches such as additive models $f(\mathbf{x}) = \sum_{k=1}^d f_k(x_k)$ or block-additive models offer a natural route to mitigate the curse of dimensionality: each component uses a one-dimensional grid, so the total cost scales as $O(dN^3)$ rather than $O(N^{3d})$. Extending Theorem~3.4 to this setting would require studying when the finite-rank additive structure preserves the parent-GP contraction rate.

\textit{Adaptive smoothness for SPDE.} Theorem~3.3 uses the oracle bandwidth $\kappa_n$ depending on the unknown smoothness $\alpha$. Whether full adaptation over both $N$ and $\kappa$ is achievable for the SPDE construction remains open. The challenge is that the sparse precision structure of the SPDE prior depends critically on the smoothness parameter $\beta$, so hierarchical modeling of $\beta$ would require understanding how the sparsity pattern changes with $\beta$. A systematic comparison with sparse variational GP methods \citep{nieman2025adaptive} in terms of both computation and statistical accuracy also merits further study.

\appendix

In the Appendix, we provide auxiliary results and proofs for the main theorems in the paper
``\papertitle''.

Appendix \ref{appendix:auxil} contains auxiliary lemmas used in the proofs.
Appendix \ref{appendix:SPDE} proves Theorem 3.3.
Appendix \ref{appendix:GPI} proves Theorem 3.4. Since the supplement introduces a considerable amount of notation, Table~\ref{tab:notation} collects the operators, index sets, and other frequently used symbols, together with where each is defined.

\begingroup
\footnotesize
\setlength{\tabcolsep}{5pt}
\renewcommand{\arraystretch}{0.95}
\begin{longtable}{@{}l p{8.6cm} l@{}}
\caption{Notation used throughout the supplement.}\label{tab:notation}\\
\toprule
Symbol & Meaning & Defined in \\
\midrule
\endfirsthead
\multicolumn{3}{@{}l}{\itshape Table \ref{tab:notation} (continued)}\\
\toprule
Symbol & Meaning & Defined in \\
\midrule
\endhead
\bottomrule
\endfoot
\multicolumn{3}{@{}l}{\emph{Indices, grids, and bases}}\\
$d$ & Input dimension (fixed) & Section~\ref{appendix:auxil} \\
$N$ & Resolution level; regular grid $\{\bj/N\}$ & Section~\ref{appendix:auxil} \\
$\cJ_N$ & Multi-index set, $|\cJ_N|=(N+1)^d$ & \eqref{eq:multiindex_set_supp} \\
$\bi,\bj,\bk$ & Multi-indices in $\cJ_N$ & Section~\ref{appendix:auxil} \\
$\psi_{j_r}$ & One-dimensional hat basis function & Section~\ref{appendix:auxil} \\
$\Psi_{\bj}$ & Tensor-product hat basis $\prod_r\psi_{j_r}$ & \eqref{eq:psi_tensor_def_supp} \\
$V_N$ & $\mathrm{span}\{\Psi_{\bi}:\bi\in\cJ_N\}$ & Section~\ref{appendix:SPDE_technical} \\
\addlinespace
\multicolumn{3}{@{}l}{\emph{Operators}}\\
$\cR_N$ & Restriction operator (nodal values) & Section~\ref{appendix:SPDE_technical} \\
$\cE_N$ & Extension operator (coefficients $\to$ piecewise-linear) & \eqref{eqn: extension operator} \\
$\cI_N=\cE_N\circ\cR_N$ & Multilinear interpolation operator & \eqref{eqn: interpolation operator} \\
$\cT_J$ & Spectral cut-off operator at level $J$ & \eqref{eq:SPDE_trunc_operator_app} \\
$S$ & Cut-off function, supported on $[-1,1]$ & Section~\ref{sec:RKHS} \\
$\cL_\kappa=\kappa^2\,\mathrm{id}-\Delta$ & SPDE operator (Neumann BC) & Section~\ref{appendix:SPDE_technical} \\
\addlinespace
\multicolumn{3}{@{}l}{\emph{Finite-element matrices}}\\
$C^{(c)},\,G$ & 1D consistent mass / stiffness matrices & Section~\ref{appendix:SPDE_technical} \\
$C^{(c)}_d,\,G^{(c)}_d$ & $d$-dimensional consistent FE matrices & Section~\ref{appendix:SPDE_technical} \\
$C$ & Mass-lumped 1D mass matrix & Section~\ref{appendix:SPDE_technical} \\
$C_d=C^{\otimes d}$ & Mass-lumped $d$-dim mass matrix & Section~\ref{appendix:SPDE_technical} \\
$G_d$ & Mass-lumped $d$-dim stiffness matrix & Section~\ref{appendix:SPDE_technical} \\
$L_{N,d},\,L^{(c)}_{N,d}$ & Discrete Galerkin operators (lumped / consistent) & Section~\ref{appendix:SPDE_technical} \\
$Q_d$ & Scale-corrected precision $\kappa^{-(2\beta-d)}C_dL_{N,d}^{\beta}$ & \eqref{eq:SPDE_Qd_app} \\
\addlinespace
\multicolumn{3}{@{}l}{\emph{Spectral and eigen-quantities}}\\
$e_{\bk}$ & Eigenbasis of $\cL_\kappa$ & Section~\ref{appendix:SPDE_technical} \\
$v_{\bk}=\cR_N e_{\bk}$ & Nodal-value vector of $e_{\bk}$ & Lemma~\ref{lem:SPDE_eigs_d_app} \\
$\mu_{\bk}$ & Generalized eigenvalue, $G_dv=\mu C_dv$ & Lemma~\ref{lem:SPDE_eigs_d_app} \\
$\lambda_{\bk}$ & Discrete precision eigenvalue, $Q_dv=\lambda C_dv$ & Lemma~\ref{lem:SPDE_eigs_d_app} \\
$\Lambda_{\bk}$ & Continuous (Mat\'ern) precision eigenvalue & Section~\ref{appendix:SPDE_technical} \\
$\gamma_{\bk}=v_{\bk}^\top C_d v_{\bk}$ & $C_d$-orthogonality constants, $1\le\gamma_{\bk}\le 2^d$ & Lemma~\ref{lem:SPDE_gamma_1d_app} \\
$z_{\bk}$ & Spectral coefficients of $f$ (or $f^*$) & Section~\ref{appendix:SPDE_technical} \\
\addlinespace
\multicolumn{3}{@{}l}{\emph{Prior, spaces, and RKHS}}\\
$w^N$ & Coefficient vector, $w^N\mid\kappa,N\sim\cN(0,Q_d^{-1})$ & \eqref{eq:SPDE_Qd_app} \\
$f_N=\cE_N w^N$ & Finite-rank prior field & Section~\ref{appendix:SPDE_technical} \\
$\beta$ & SPDE order parameter & \eqref{eq:SPDE_Qd_app} \\
$\kappa\;(\kappa_n)$ & Inverse-bandwidth / scale parameter & Section~\ref{appendix:SPDE Main proof} \\
$\bbH(\kappa,N),\,\bbH_1(\kappa,N)$ & RKHS and its unit ball & Section~\ref{appendix:SPDE_technical} \\
$A_R^\infty$ & Sobolev-type ellipsoid of radius $R$ & \eqref{eq:ARk_def_supp} \\
$C^\alpha,\;\Besova$ & H\"older / Besov smoothness classes & Section~\ref{sec:RKHS} \\
$f^*,\;C_0=\|f^*\|_{\Besova}$ & Target function and its Besov norm & Lemma~\ref{lem:SPDE_approx_app} \\
\addlinespace
\multicolumn{3}{@{}l}{\emph{Sieve and contraction (Appendix \ref{appendix:SPDE})}}\\
$\epsilon_n$ & Contraction rate, $n^{-\alpha/(2\alpha+d)}$ up to logs & Section~\ref{appendix:SPDE Main proof} \\
$\cB_n,\,\cB_n^{(0)},\,\cB_{n,N}$ & Sieve and its components & \eqref{eq:SPDE_oracle_sieve_app} \\
$M_n,\,R_n,\,N_n$ & Sieve scale / radius / resolution cutoff & Section~\ref{appendix:SPDE Main proof} \\
$\bbB_1$ & Unit ball of $C[0,1]^d$ in $\|\cdot\|_\infty$ & Section~\ref{appendix:SPDE Main proof} \\
\addlinespace
\multicolumn{3}{@{}l}{\emph{GPI approach (Appendix \ref{appendix:GPI})}}\\
$g$ & Squared-exponential parent GP & \eqref{eqn: RBF kernel function} \\
$d_g$ & Canonical pseudometric of $g$ & Section~\ref{appendix:GPI_technical} \\
$B^\infty_\epsilon,\,B^g_\epsilon$ & $\ell_\infty$- and $d_g$-balls of pairs $(s,t)$ & \eqref{eq:B_sets_supp} \\
$C_{\mathrm{osc}},C_{\mathrm{vol}},c_{\mathrm{Bor}}$ & Constants in modulus-of-continuity / Borell bounds & Section~\ref{appendix:GPI_technical} \\
$\cE(\epsilon)$ & Modulus-of-continuity event & Section~\ref{appendix:gpi_sieve_universalN} \\
$\cG(\epsilon),\,\widetilde{\cG}(\epsilon),\,\cF^*(\epsilon)$ & Base sieve, stabilized sieve, final sieve & \eqref{eq:gpi_stabilizedG_def} \\
$\kappa_\epsilon,\,N_0(\epsilon),\,N_1(\epsilon)$ & Scale / resolution cutoffs & Section~\ref{appendix:GPI_main} \\
\end{longtable}
\endgroup

\section{Auxiliary results}\label{appendix:auxil}

Throughout the supplement, the dimension $d$ is assumed fixed. For an integer $N\ge 1$, recall the multi-index set
\begin{equation}\label{eq:multiindex_set_supp}
\cJ_N:=\{(j_1,\ldots,j_d): j_r\in\{0,1,\ldots,N\} \text{ for } r \in [d]\}.
\end{equation}
Note that $|\cJ_N| = (N+1)^d$. We write bold indices such as $\bi,\bj,\bk\in\cJ_N$. Also, as in the main document, $\{\psi_{j_r}\}_{j_r=0}^N$ throughout this document denote one-dimensional hat basis functions given by
\begin{equation}
    \psi_{j_r}(x) = \Big(1 - N\Big|x - \frac{j_r}{N}\Big|\Big)^+,
\end{equation}
and $d$-dimensional hat basis is defined as a tensor product:
\begin{equation}\label{eq:psi_tensor_def_supp}
\Psi_{\bj}(x)=\prod_{r=1}^d \psi_{j_r}(x_r), \qquad \bj=(j_1,\ldots,j_d)\in\cJ_N.
\end{equation}

\subsection{Properties of finite basis expansion}\label{appendix:general}

We begin with the proof of Lemma 3.1, which describes the RKHS norm associated with $f_N$.

\begin{proof}[Proof of Lemma 3.1]
We use Theorem 4.1 of \citet{van2008reproducing}, which characterizes the RKHS associated with a (finite) Gaussian series representation.

Define the feature vector $\Psi(x)=(\psi_j(x))_{j\in\cJ_N}\in\bbR^{(N+1)^d}$, so that
$f_N(x)=\sum_{j\in\cJ_N} w_j \psi_j(x) = \Psi(x)^\top w^N$.
Then the covariance kernel of $f_N$ is
\begin{equation}\label{eq:cov_kernel_series_supp}
K(x_1,x_2)=\bbE\{f_N(x_1)f_N(x_2)\}=\Psi(x_1)^\top \Sigma_{N, \theta} \Psi(x_2).
\end{equation}

Let $\Sigma_{N, \theta} = V\Lambda V^\top$ be an eigendecomposition, with $V$ orthogonal and
$\Lambda=\mathrm{diag}(\lambda_1,\ldots,\lambda_m)$ where $m=(N+1)^d$ and $\lambda_\ell>0$ for $\ell \in [m]$.
Set $\widetilde\Psi(x)=V^\top\Psi(x)$ and write
$\widetilde\Psi(x)=(\widetilde\psi_1(x),\ldots,\widetilde\psi_m(x))^\top$.
Then
\begin{equation}\label{eq:cov_kernel_eig_supp}
K(x_1,x_2)=\widetilde\Psi(x_1)^\top \Lambda \widetilde\Psi(x_2)
=\sum_{\ell=1}^m \lambda_\ell\,\widetilde\psi_\ell(x_1)\widetilde\psi_\ell(x_2).
\end{equation}
By Theorem 4.1 of \citet{van2008reproducing}, the RKHS $\bbH$ consists of all functions
$h(x)=\sum_{\ell=1}^m a_\ell \widetilde\psi_\ell(x)$ equipped with norm
\begin{equation}\label{eq:rkhs_norm_eig_supp}
\|h\|_{\bbH}^2=\sum_{\ell=1}^m a_\ell^2/\lambda_\ell.
\end{equation}
Since $\{\widetilde\psi_\ell\}$ span the same space as $\{\psi_j\}$, we have
$\bbH=\mathrm{span}\{\psi_j:j\in\cJ_N\}$.

Now let $h=\sum_{j\in \cJ_N}v_j\psi_j$ and write $v^N=Va$ for some $a \in \bbR^m$. Then
\begin{equation*}
\|h\|_{\bbH}^2 = a^\top\Lambda^{-1}a =
(v^N)^\top V\Lambda^{-1}V^\top v^N =
(v^N)^\top\Sigma_{N,\theta}^{-1}v^N .
\end{equation*}
This concludes the proof.
\end{proof}

We next prove Lemma 3.2, which gives a sufficient condition on $N$ for approximating $f^*$ by a tensor-product hat basis.

\begin{proof}[Proof of Lemma 3.2]
With $d$-dimensional hat basis functions defined in \eqref{eq:psi_tensor_def_supp}, define the multilinear interpolation of $f^*$ on the regular grid as
\begin{equation}\label{eq:interp_def_general_d}
    f_N(x)=\sum_{\bj\in\cJ_N} f^*(\bj/N)\,\Psi_{\bj}(x),
    \qquad x\in[0,1]^d.
\end{equation}
By construction, $f_N\in\mathrm{span}\{\Psi_\bj:\bj\in\cJ_N\}$ and $f_N(\bj/N)=f^*(\bj/N)$ for all grid points.

Write $h=1/N$ and fix a cell $Q=\prod_{r=1}^d [i_r h,(i_r+1)h]$.
For any $x\in Q$, $f_N(x)$ is a convex combination of the vertex values $\{f^*(u):u\in\mathcal{V}(Q)\}$, i.e.
$f_N(x)=\sum_{u\in\mathcal{V}(Q)} \theta_u(x) f^*(u)$ with $\theta_u(x)\ge 0$ and $\sum_u\theta_u(x)=1$,
and every vertex $u\in\mathcal{V}(Q)$ satisfies $\|x-u\|_\infty\le h$.

We consider three regimes.

\medskip\noindent
\textbf{Case 1: $0<\alpha\le 1$.}
Since $f^*\in C^\alpha([0,1]^d)$, there exists $C_{f^*}>0$ such that
$|f^*(x)-f^*(u)|\le C_{f^*}\|x-u\|_\infty^\alpha \le C_{f^*}h^\alpha$ for all $u\in\mathcal{V}(Q)$.
Therefore,
\begin{equation}\label{eq:interp_err_case1_supp}
|f^*(x)-f_N(x)|
=\Big|\sum_{u\in\mathcal{V}(Q)}\theta_u(x)\{f^*(x)-f^*(u)\}\Big|
\le \sum_u\theta_u(x)\,C_{f^*}h^\alpha
= C_{f^*}h^\alpha.
\end{equation}

\medskip\noindent
\textbf{Case 2: $1<\alpha\le 2$.}
Then $\nabla f^*$ exists and is $(\alpha-1)$-H\"older. For each vertex $u\in\mathcal{V}(Q)$,
a first-order Taylor expansion around $x$ yields
\begin{equation}\label{eq:taylor_1st_supp}
f^*(u)=f^*(x)+\nabla f^*(x)^\top(u-x)+R_u,
\qquad |R_u|\le C_{f^*}\|u-x\|_\infty^\alpha\le C_{f^*}h^\alpha.
\end{equation}
Multilinear interpolation reproduces affine functions exactly, hence
$\sum_u\theta_u(x)\{f^*(x)+\nabla f^*(x)^\top(u-x)\}=f^*(x)$.
Thus
\begin{equation}\label{eq:interp_err_case2_supp}
|f^*(x)-f_N(x)|
=\Big|\sum_{u\in\mathcal{V}(Q)}\theta_u(x)R_u\Big|
\le \sum_u\theta_u(x)\,C_{f^*}h^\alpha
= C_{f^*}h^\alpha.
\end{equation}

\medskip\noindent
\textbf{Case 3: $\alpha>2$.}
Then all second-order partial derivatives (including mixed partial derivatives) exist and are bounded.
A second-order Taylor expansion gives
\begin{equation}\label{eq:taylor_2nd_supp}
f^*(u)=f^*(x)+\nabla f^*(x)^\top(u-x)+\tfrac12 (u-x)^\top Hf^*(\xi_u)(u-x),
\end{equation}
for some $\xi_u$ on the line segment between $x$ and $u$.
Again the affine part is reproduced exactly by multilinear interpolation, hence
\begin{equation}\label{eq:interp_err_case3_supp}
|f^*(x)-f_N(x)|
\le \frac12 \sup_{z\in[0,1]^d}\|Hf^*(z)\|_{\mathrm{op}}
\sum_{u\in\mathcal{V}(Q)}\theta_u(x)\|u-x\|_2^2
\lesssim h^2.
\end{equation}

Taking the supremum over $x\in[0,1]^d$ yields
\begin{equation}\label{eq:interp_err_final_supp}
\|f_N-f^*\|_\infty \le C\,h^{\alpha\land 2}=C\,N^{-(\alpha\land 2)}
\end{equation}
for a constant $C>0$ depending only on $f^*$ (through its H\"older norm) and $d$.
Therefore, choosing $N\ge (C/\epsilon)^{1/(\alpha\land 2)}$ ensures $\|f_N-f^*\|_\infty<\epsilon$.
\end{proof}


\section{Proof of Theorem 3.3}\label{appendix:SPDE}

The proof of Theorem 3.3 compares the finite-element SPDE prior with the corresponding continuous Mat\'{e}rn prior under the same choice of $\kappa_n$. The comparison is mainly carried out between the eigenvalues of the differential operators in the SPDE defining the parent GP and the generalized eigenvalues of the finite-rank GP. Lemma \ref{lem:SPDE_eigs_1d_app} to Lemma \ref{lem:SPDE_eigs_d_app} establish the comparison between the discrete precision eigenvalues and the continuous precision eigenvalues. Lemma \ref{lem:SPDE_entropy_H1_app} to Lemma \ref{lem:SPDE_RKHS_norm_app} use this comparison to control the RKHS approximation term, centered small-ball probability, and entropy of the finite-rank prior. Section B.1 collects these technical lemmas, and Section B.2 verifies the three posterior-contraction conditions.

\subsection{Technical lemmas}\label{appendix:SPDE_technical}

This section gives the finite-element matrix details used in the compact SPDE construction in Section 2.1 of the main manuscript.

For the SPDE method, we use several representations of the same function. For $f \in \Besova$, the spectral representation is written as
\begin{equation*}
f = \sum_k z_k e_k,
\end{equation*}
where $\{e_k\}$ denotes the eigenbasis functions of $\cL_\kappa = \kappa^2 \id - \Delta$ under Neumann boundary condition.

The {\it restriction operator} $\cR_N$ maps a continuous function to its vector of nodal values,
\begin{equation}\label{eqn: restriction operator}
(R_N f)_i = f\biggl({\frac{i_1}{N}}, \ldots, \frac{i_d}{N}\biggr),
\qquad i \in \cJ_N.
\end{equation}
Next, the {\it extension operator} $\cE_N$ maps the nodal values into a piecewise linear function:
\begin{equation}\label{eqn: extension operator}
f_N = \cE_N w = \sum_{\bi \in \cJ_N} w_i \Psi_i,
\end{equation}
where $\{\Psi_\bi : \bi \in \cJ_N\}$ is the tensor-product hat basis defined in \eqref{eq:psi_tensor_def_supp}. 
Finally, the {\it multilinear interpolation operator} on the regular grid is defined as
\begin{equation}\label{eqn: interpolation operator}
    \cI_N = \cE_N \circ \cR_N.
\end{equation}
The interpolation operator $\cI_N$ maps a continuous function into a piecewise linear function sharing the same nodal values on the grid.

Let
\begin{equation*}
V_N=\operatorname{span}\{\Psi_{\bi}:\bi\in\cJ_N\}.
\end{equation*}

Recall that the one-dimensional consistent mass matrix and stiffness matrix are defined as
\begin{equation}\label{eq:CG_consistent_def_supp}
C^{(c)}_{ij}=\int_0^1 \psi_i(t)\psi_j(t)dt,
\qquad
G_{ij}=\int_0^1 \psi_i'(t)\psi_j'(t)dt,
\qquad 0\le i,j\le N.
\end{equation}
The corresponding $d$-dimensional consistent finite-element matrices are
\begin{equation}\label{eq:CdGd_def_supp}
(C_d^{(c)})_{\bi,\bj}
=
\int_{\Omega}\Psi_{\bi}(\bx)\Psi_{\bj}(\bx)d\bx,
\qquad
(G_d^{(c)})_{\bi,\bj}
=
\int_{\Omega}\nabla\Psi_{\bi}(\bx)^\top\nabla\Psi_{\bj}(\bx)d\bx .
\end{equation}
For coefficient vectors $u=(u_{\bi})_{\bi\in\cJ_N}$ and $v=(v_{\bi})_{\bi\in\cJ_N}$, write
\begin{equation*}
u_N=\sum_{\bi\in\cJ_N}u_{\bi}\Psi_{\bi},
\qquad
v_N=\sum_{\bi\in\cJ_N}v_{\bi}\Psi_{\bi}.
\end{equation*}Then
\begin{equation*}
v^\top C_d^{(c)}u=\langle v_N,u_N\rangle_{L^2(\Omega)},
\qquad
v^\top G_d^{(c)}u
=
\int_\Omega \nabla v_N(\bx)^\top\nabla u_N(\bx)d\bx .
\end{equation*}

More precisely, for $u_N=\sum_{\bi\in\cJ_N}u_{\bi}\Psi_{\bi}$, define
$z_N=\sum_{\bi\in\cJ_N}z_{\bi}\Psi_{\bi}\in V_N$ by
\begin{equation*}
\langle v_N,z_N\rangle_{L^2(\Omega)}
=
\int_\Omega \nabla v_N(\bx)^\top\nabla u_N(\bx)d\bx
\qquad \mbox{for all } v_N\in V_N.
\end{equation*}
In coefficient coordinates this identity is
\begin{equation*}
v^\top C_d^{(c)}z=v^\top G_d^{(c)}u
\qquad \mbox{for all } v\in\mathbb{R}^{(N+1)^d}.
\end{equation*}
Hence
\begin{equation*}
C_d^{(c)}z=G_d^{(c)}u,
\qquad
z=\big(C_d^{(c)}\big)^{-1}G_d^{(c)}u,
\end{equation*}
because $C_d^{(c)}$ is positive definite under linearly independent $\{\Psi_\bi\}$. Therefore
$\big(C_d^{(c)}\big)^{-1}G_d^{(c)}$ is the consistent Galerkin matrix for
$-\Delta$ in coefficient coordinates, and
\begin{equation*}
L_{N, d}^{(c)} = \big(C_d^{(c)}\big)^{-1}\big(\kappa^2 C_d^{(c)}+G_d^{(c)}\big)
=
\kappa^2 I+\big(C_d^{(c)}\big)^{-1}G_d^{(c)}
\end{equation*}
is the matrix leading to coefficient-coordinate version of $\cL_\kappa=\kappa^2\id-\Delta$.

By the tensor-product form of $\Psi_{\bi}$,
\begin{equation}\label{eq:Cd_consistent_kron_supp}
C_d^{(c)}=(C^{(c)})^{\otimes d},
\end{equation}
and
\begin{equation}\label{eq:Gd_consistent_kron_supp}
G_d^{(c)}
=
\sum_{r=1}^d
\Big((C^{(c)})^{\otimes(r-1)}\otimes G\otimes (C^{(c)})^{\otimes(d-r)}\Big).
\end{equation}
In practice, the SPDE prior uses mass lumping for $C^{(c)}$ to induce sparsity in the final precision matrix. In one dimension, only the mass matrix is replaced by
\begin{equation*}
C := \diag(C^{(c)}\mathbf{1}_{N+1})
=
\frac{1}{2N}\diag(1,2,\ldots,2,1),
\end{equation*}
while the one-dimensional stiffness matrix $G$ is unchanged. For the tensor-product basis, the mass-lumped matrix is
\begin{equation*}
C_d
:=
\diag(C_d^{(c)}\mathbf{1}_{(N+1)^d})
=
C^{\otimes d},
\end{equation*}
and the stiffness matrix uses $C^{(c)}$ replaced with $C$:
\begin{equation}\label{eq:Gd_lumped_kron_supp}
G_d
:=
\sum_{r=1}^d
\Big(C^{\otimes(r-1)}\otimes G\otimes C^{\otimes(d-r)}\Big).
\end{equation}
Thus $G_d=G_d^{(c)}$ when $d=1$, while $G_d$ generally differs from $G_d^{(c)}$ when $d\ge2$. $L_{N, d}$ is obtained from $L_{N, d}^{(c)}$ by replacing $C_d^{(c)}$ and $G_d^{(c)}$ to $C_d$ and $G_d$, respectively.
From this point onward, $C_d$ and $G_d$ denote the mass-lumped matrices used in the SPDE prior. The main manuscript suppresses the subscript $d$ and writes these matrices as $C$ and $G$. 

Let $\beta\in\bbN$. The scale-corrected precision matrix is
\begin{equation}\label{eq:SPDE_Qd_app}
    Q_d=
    \kappa^{-(2\beta-d)}C_dL_{N, d}^\beta
\end{equation}
and the conditional prior is
\begin{equation}\label{eq:SPDE_prior_app}
    w^N\mid \kappa,N\sim \cN(0,Q_d^{-1}),
    \qquad
    f_N=\cE_N w^N.
\end{equation}

Next, we record the continuous spectral coefficients corresponding to the parent Mat\'ern
SPDE. Since
\begin{equation*}
\cL_\kappa e_k
=
\left(\kappa^2+\pi^2\|k\|_2^2\right)e_k,
\qquad
\cL_\kappa=\kappa^2\mathrm{id}-\Delta,
\end{equation*}
the SPDE
\begin{equation*}
\kappa^{-(\beta-d/2)}\cL_\kappa^{\beta/2}f = W
\end{equation*}
implies that, for $f=\sum_k z_k e_k$,
\begin{equation*}
z_k \sim N(0,\Lambda_k^{-1}),
\end{equation*}
where
\begin{equation*}
\Lambda_k
=
\kappa^{-(2\beta-d)}
\left(
\kappa^2+\pi^2\|k\|_2^2
\right)^\beta .
\end{equation*}
Equivalently,
\begin{equation*}
\kappa^{-(2\beta-d)}\cL_\kappa^\beta e_k = \Lambda_k e_k .
\end{equation*}

We then introduce the results related to eigenvectors $e_{\bk}$. We start with $d=1$. For $k \ge 1$, since
\begin{equation*}
    e_k(s) - e_k(t) = -2\sqrt{2} \sin(k\pi(s+t)/2) \sin(k\pi(s-t)/2),
\end{equation*}
combined with $|\sin x| \le \min(|x|, 1)$, we obtain $|e_k(s) - e_k(t)| \le \min(2\sqrt{2}, \sqrt{2} k\pi|s-t|)$ for $0 \le k \le N$. This can be extended to general dimension as follows.
\begin{lemma}\label{lem:e k Holder continuity}
    For $\bs, \bt \in [0, 1]^d$ and $0 < \rho < 1$,
    \begin{equation*}
        |e_{\bk}(\bs) - e_{\bk}(\bt)| \le d \cdot  2^{d/2 + 1 - \rho} \|\pi \bk\|^{\rho}_2 \| \bs - \bt\|_{\infty}^{\rho}.
    \end{equation*}
\end{lemma}
\begin{proof}
    First, observe 
    \begin{align*}
        |e_{\bk}(\bs) - e_{\bk}(\bt)| &\le \sum_{r=1}^d |e_{\bk}(s_1, \ldots, s_r, t_{r+1}, \ldots, t_d) - e_{\bk}(s_1, \ldots, s_{r-1}, t_{r}, \ldots, t_d)|\\
        & \le (\sqrt{2})^{d-1}\sum_{r=1}^d |e_{k_r}(s_r) - e_{k_r}(t_r)|.
    \end{align*}
    Since $\min(a, b) \le a^{1-\rho} b^{\rho}$ holds for all $0 < \rho < 1$ and $a, b \ge 0$, we observe
    \begin{equation*}
        |e_k(s) - e_k(t)| \le \min(2\sqrt{2}, \sqrt{2} k\pi|s-t|) \le  (2\sqrt{2})^{1-\rho} (\sqrt{2} k\pi |s-t|)^{\rho}.
    \end{equation*}
    Combining the above two display, we obtain
    \begin{equation*}
        |e_{\bk}(\bs) - e_{\bk}(\bt)| \le d \cdot 2^{d/2 + 1 - \rho} \|\pi \bk \|_2^\rho \|\bs - \bt \|_{\infty}^{\rho}.
    \end{equation*}
\end{proof}
Additionally, we define the following ellipsoid :
\begin{equation}\label{eq:ARk_def_supp}
    A_{R}^\infty
    =
    \Big\{\sum_{\bk\in\bbN_0^d} z_{\bk} e_{\bk}:\ \sum_{\bk\in\bbN_0^d} (1+\|\pi\bk\|_2^2)^\beta z_{\bk}^2 \le R^2\Big\}.
\end{equation}
Then the following Lemma guarantees the H\"{o}lder continuity of the elements in $A_R^{\infty}$.
\begin{lemma}\label{lem:holder continuity of AR}
    For $g \in A_R^{\infty}$, for $0 < \rho < \min(\beta - d/2, 1)$,
    \begin{equation}\label{eq:holder continuity of AR}
        |g(\bs) - g(\bt)| \le C_{d, \beta, \rho} R\|s-t\|^{\rho}_{\infty}
    \end{equation}
    holds for all $\bs, \bt \in [0, 1]^d$, with a constant $C_{d, \beta, \rho} >0$ only depending on $d, \beta$, and $\rho$.
\end{lemma}
\begin{proof}
    From Lemma \ref{lem:e k Holder continuity}, for all $0 < \rho < 1$,
    \begin{align*}
        |g(\bs) - g(\bt)| \le \sum_{\bk} |z_{\bk}||e_{\bk}(\bs) - e_{\bk}(\bt)| \le d \cdot 2^{d/2 + 1 - \rho} \|\bs - \bt\|_{\infty}^{\rho}\sum_{\bk} |z_{\bk}| \|\pi \bk \|_2^{\rho} 
    \end{align*}
    holds. From Cauchy--Schwarz inequality, we obtain
    \begin{equation*}
        \sum_{\bk} |z_{\bk}| \|\pi \bk \|_2^{\rho}  \le \Big( \sum_{\bk} (1+\|\pi\bk\|_2^2)^\beta z_{\bk}^2  \Big)^{1/2} \Bigg(\sum_{\bk} \frac{\| \pi \bk\|_2^{2\rho}}{(1 + \|\pi \bk \|_2^2)^\beta} \Bigg)^{1/2}.
    \end{equation*}
    For $\rho$ satisfying $\rho < \beta - d/2$, the last term converges to a constant only depending on $\rho, \beta$, and $d$. Therefore, for 
    \begin{equation*}
        C_{d, \beta, \rho} = \Bigg(\sum_{\bk} \frac{\| \pi \bk\|_2^{2\rho}}{(1 + \|\pi \bk \|_2^2)^\beta} \Bigg)^{1/2},
    \end{equation*}
    \eqref{eq:holder continuity of AR} holds.
\end{proof}

Subsequently, we introduce the results on finite rank SPDE representation. First, write $v_{\bk}= \cR_N e_{\bk}$ as the nodal value vector of $e_{\bk}$ with $\cR_N$ defined in \eqref{eqn: restriction operator}. The following lemma provides the link among the associated eigenvalues, eigenbasis, and eigenvectors.

\begin{lemma}[Generalized eigenpairs in $d=1$]\label{lem:SPDE_eigs_1d_app}
Let $d=1$ and $v_k=\cR_N e_k\in\bbR^{N+1}$ for $0\le k\le N$. Then
\begin{equation}\label{eq:SPDE_GC_eig_1d_app}
    G v_k = \mu_k\, C v_k,
    \qquad
    \mu_k = 4N^2\sin^2\Big(\frac{k\pi}{2N}\Big),
\end{equation}
and if $Q$ denotes $Q_d$ in \eqref{eq:SPDE_Qd_app} with $d=1$, then
\begin{equation}\label{eq:SPDE_QC_eig_1d_app}
    Q v_k = \lambda_k\, C v_k,
    \qquad
    \lambda_k
    =
    \kappa^{-(2\beta-1)}
    \Big\{\kappa^2+\mu_k\Big\}^{\beta}.
\end{equation}
Moreover, we have
\begin{equation}\label{eq:SPDE_lambda_compare_1d_app}
    \lambda_k\le \Lambda_k \le \lambda_k(\pi/2)^{2\beta}.
\end{equation}
\end{lemma}

\begin{proof}
First, we observe $Gv_0 = 0_{N+1} = \mu_0 Cv_0$. For $k \ge 1$, let $(v_k)_i = e_k(i/N)$ for $0 \le i \le N$ and $\theta_k = k\pi / N$. Observe
\begin{align}
    (Gv_k)_0 &= N((v_k)_0 - (v_k)_1), \\
    (Gv_k)_i &= N(-(v_k)_{i-1} +2(v_k)_i - (v_k)_{i+1}), \quad 1\le i\le N-1,\\ 
    (Gv_k)_N &= N((v_k)_N - (v_k)_{N-1}).
\end{align}
For $1\le i\le N-1$, from $\cos(a-b) + \cos(a+b) = 2\cos a \cos b$, we obtain
\begin{equation*}
    (Gv_k)_i = 2N(1 - \cos \theta_k) (v_k)_i = 4N \sin^2(\theta_k/2) (v_k)_i = 4N^2 \sin^2(\theta_k /2) (Cv_k)_i.
\end{equation*}
For $i = 0$,
\begin{equation*}
    (Gv_k)_0=  \sqrt{2}N(1 - \cos \theta_k) = 2\sqrt{2}N\sin^2(\theta_k/2) = 4N^2\sin^2(\theta_k/2) (Cv_k)_0.
\end{equation*}
For $i = N$,
\begin{equation*}
    (Gv_k)_N=  \sqrt{2}N(-1)^k(1 - \cos \theta_k) = 2\sqrt{2}N(-1)^k\sin^2(\theta_k/2) = 4N^2\sin^2(\theta_k/2) (Cv_k)_N.
\end{equation*}
Therefore, we obtain $Gv_k = \mu_k Cv_k$.
For \eqref{eq:SPDE_QC_eig_1d_app}, let $A=\kappa^2 C+G$ (the $d=1$ case). From \eqref{eq:SPDE_GC_eig_1d_app} we have
$A v_k = (\kappa^2+\mu_k)\,C v_k$, hence $(C^{-1}A)v_k=(\kappa^2+\mu_k)v_k$ and therefore
$C(C^{-1}A)^{\beta}v_k = (\kappa^2+\mu_k)^{\beta}\,C v_k$.
Using \eqref{eq:SPDE_Qd_app} with $d=1$ gives
$Q v_k=\kappa^{-(2\beta-1)}\,C(C^{-1}A)^{\beta}v_k
=\kappa^{-(2\beta-1)}(\kappa^2+\mu_k)^{\beta}\,C v_k$,
which is \eqref{eq:SPDE_QC_eig_1d_app}. Lastly, \eqref{eq:SPDE_lambda_compare_1d_app} is directly obtained from $(2x)/\pi\le \sin(x)\le x$ on $[0,\pi/2]$ with $x=k\pi/(2N)$.
\end{proof}

\begin{lemma}[$C$-orthogonality and the constants $\gamma_k$ in $d=1$]\label{lem:SPDE_gamma_1d_app}
Let $d=1$ and $v_k=\cR_N e_k$. Then
\begin{equation}\label{eq:SPDE_Corth_1d_app}
    v_k^\top C v_\ell = \gamma_k\,\ind_{\{k=\ell\}},
\end{equation}
where
\begin{equation}\label{eq:SPDE_gamma_def_1d_app}
    \gamma_k=
    \begin{cases}
        1, & 0\le k < N,\\
        2, & k=N.
    \end{cases}
\end{equation}
\end{lemma}

\begin{proof}
Recall that
\begin{equation}\label{eq:C_lumped_explicit_supp}
C=\frac{1}{2N}\,\mathrm{diag}(1,2,\ldots,2,1).
\end{equation}
Hence for vectors $a,b\in\bbR^{N+1}$,
\begin{equation}\label{eq:disc_innerprod_supp}
a^\top C b=\frac{1}{2N}\Big\{a_0 b_0+a_N b_N+2\sum_{i=1}^{N-1} a_i b_i\Big\}
=\frac{1}{N}\sum_{i=0}^N w_i\,a_i b_i,
\end{equation}
where $w_0=w_N=1/2$ and $w_i=1$ for $1\le i\le N-1$.

Let $\theta_i=\pi i/N$. For $k\ge 1$, $v_k(i)=\sqrt{2}\cos(k\theta_i)$ and $v_0(i)=1$.
Define $T_m=\sum_{i=0}^N w_i\cos(m\theta_i)$ for integers $m\ge 0$.
A direct geometric-series calculation gives
\begin{equation}\label{eq:cos_sum_identity_supp}
T_0=\sum_{i=0}^N w_i = N,\qquad
T_m=0\ \ \text{for }m=1,\ldots,2N-1,\qquad
T_{2N}=N.
\end{equation}

If $k=0$ and $\ell\ge 1$, then $v_0^\top C v_\ell=(\sqrt{2}/N)T_\ell=0$ by \eqref{eq:cos_sum_identity_supp}.
Now take $k,\ell\ge 1$. Using $\cos a\cos b=\tfrac12\{\cos(a-b)+\cos(a+b)\}$ and \eqref{eq:disc_innerprod_supp},
\begin{align}
v_k^\top C v_\ell
&=\frac{2}{N}\sum_{i=0}^N w_i\cos(k\theta_i)\cos(\ell\theta_i)
\\&=\frac{1}{N}\sum_{i=0}^N w_i\Big\{\cos((k-\ell)\theta_i)+\cos((k+\ell)\theta_i)\Big\}
=\frac{1}{N}\{T_{|k-\ell|}+T_{k+\ell}\}.\label{eq:coscos_reduce_supp}
\end{align}
If $k\neq \ell$, then $|k-\ell|\in\{1,\ldots,N\}$ and $k+\ell\in\{1,\ldots,2N-1\}$, so both terms vanish by \eqref{eq:cos_sum_identity_supp}.
If $k=\ell<N$, then $T_{|k-\ell|}=T_0=N$ and $T_{k+\ell}=T_{2k}=0$, hence $v_k^\top C v_k=1$.
If $k=\ell=N$, then $T_{|k-\ell|}=T_0=N$ and $T_{k+\ell}=T_{2N}=N$, hence $v_N^\top C v_N=2$.
This proves \eqref{eq:SPDE_Corth_1d_app} and \eqref{eq:SPDE_gamma_def_1d_app}.
\end{proof}

\begin{lemma}[Generalized eigenpairs of $Q_d$ and $C_d$]\label{lem:SPDE_eigs_d_app}
Let $d\ge 1$ and $v_{\bk}=\cR_N e_{\bk}$ for $\bk\in\cJ_N$. Then
\begin{equation}\label{eq:SPDE_Qd_eig_app}
    Q_d v_{\bk}=\lambda_{\bk}\, C_d v_{\bk},
    \qquad
    \lambda_{\bk}
    =
    \kappa^{-(2\beta-d)}
    \Bigg\{\kappa^2+\sum_{r=1}^d 4N^2\sin^2\Big(\frac{k_r\pi}{2N}\Big)\Bigg\}^{\beta}.
\end{equation}
Moreover, 
we have
\begin{equation}\label{eq:SPDE_lambda_compare_d_app}
    \lambda_{\bk}\le \Lambda_{\bk}\le \lambda_{\bk}(\pi/2)^{2\beta}.
\end{equation}
\end{lemma}

\begin{proof}
The identity follows from the Kronecker product structure: $v_{\bk}=v_{k_1}\otimes\cdots\otimes v_{k_d}$, $C_d=C^{\otimes d}$, and $G_d$ is a Kronecker sum in \eqref{eq:Gd_lumped_kron_supp}.
Applying Lemma \ref{lem:SPDE_eigs_1d_app} coordinatewise yields \eqref{eq:SPDE_Qd_eig_app}.
The comparison \eqref{eq:SPDE_lambda_compare_d_app} follows by applying $(2x)/\pi\le \sin(x)\le x$ for all $ \in [0, \pi/2]$ and summing over $r$.
\end{proof}

\begin{corollary}\label{cor:SPDE_RKHS_ellipsoid_app}
Let $\bbH(\kappa,N)$ be the RKHS associated with $\cE_Nw^N$ under $w^N\sim \cN(0,Q_d^{-1})$ and extension operator $\cE_N$ defined in \eqref{eqn: extension operator}. Then the unit RKHS ball is defined as
\begin{equation}\label{eq:SPDE_H1_def_app}
    \bbH_1(\kappa,N)=\{\cE_N w:\ w^\top Q_d w\le 1\}.
\end{equation}
Additionally, for any $w\in\bbR^{(N+1)^d}$, there exist coefficients $\{z_{\bk}\}$ such that $w=\sum_{\bk\in\cJ_N} z_{\bk} v_{\bk}$ and
\begin{equation}\label{eq:SPDE_ellipsoid_constraint_app}
    w^\top Q_d w
    =
    \sum_{\bk\in\cJ_N} \lambda_{\bk}\,\gamma_{\bk}\, z_{\bk}^2,
    \qquad
    \gamma_{\bk}:=v_{\bk}^\top C_d v_{\bk}.
\end{equation}
Moreover, $1 \le \gamma_{\bk} \le 2^d$ holds for all $\bk$ and $N$.
\end{corollary}

\begin{proof}
The expansion follows from $C_d$-orthogonality of $\{v_{\bk}\}$, which is the tensorization of Lemma \ref{lem:SPDE_gamma_1d_app}. The identity \eqref{eq:SPDE_ellipsoid_constraint_app} follows by combining Lemma \ref{lem:SPDE_eigs_d_app} with $v_{\bk}^\top C_d v_{\bell}=0$ for $\bk\neq\bell$, and $\gamma_{\bk} = \prod_{r=1}^d \gamma_{k_r}$.
\end{proof}

\begin{lemma}[Metric entropy of $\bbH_1(\kappa,N)$ in $\|\cdot\|_\infty$]\label{lem:SPDE_entropy_H1_app}
Assume $\kappa\ge 1$. There exists $C>0$ (depending only on $d$ and $\beta$) such that for all $\epsilon\in(0,1)$ and all $N\ge 1$,
\begin{equation}\label{eq:SPDE_entropy_bound_app}
    \log \cN\big(\epsilon,\bbH_1(\kappa,N),\|\cdot\|_\infty\big)
    \le
    C(\kappa^{\beta-d/2}/\epsilon)^{d/\beta}.
\end{equation}
\end{lemma}

\begin{proof}
Let $h=\cE_N w\in\bbH_1(\kappa,N)$ and write $w=\sum_{\bk\in\cJ_N} z_{\bk} v_{\bk}$. By \eqref{eq:SPDE_ellipsoid_constraint_app} and \eqref{eq:SPDE_lambda_compare_d_app},
\begin{equation*}
    \sum_{\bk}\Lambda_{\bk} z_{\bk}^2 \le (\pi/2)^{2\beta}\sum_{\bk}\lambda_{\bk} \gamma_{\bk} z_{\bk}^2 \le (\pi/2)^{2\beta}.
\end{equation*}
Under $\kappa\ge 1$,
\begin{equation}\label{eq:Lambda_lower_supp}
    \Lambda_{\bk}
    =
    \kappa^{-(2\beta-d)}\Big(\kappa^2+\|\pi\bk\|_2^2\Big)^\beta
    \ge
    \kappa^{-(2\beta-d)}\Big(1+\|\pi\bk\|_2^2\Big)^\beta.
\end{equation}
Hence, 
\begin{equation*}
    \sum_{\bk}(1+\|\pi\bk\|_2^2)^\beta z_{\bk}^2 \le (\pi/2)^{2\beta}\kappa^{2\beta-d}.
\end{equation*}
This leads to $\widetilde h(\bx)=\sum_{\bk} z_{\bk} e_{\bk}(\bx) \in A_{R_\kappa}^{\infty}$ with $R_\kappa=(\pi/2)^{\beta}\,\kappa^{\beta-d/2}$.
Since $\cI_N$ is a contraction in $\|\cdot\|_\infty$ (the hat basis is nonnegative and sums to one),
\begin{equation}\label{eq:entropy_reduction_supp}
    \cN\big(\epsilon,\bbH_1(\kappa,N),\|\cdot\|_\infty\big)
    \le
    \cN\big(\epsilon, A_{R_\kappa}^\infty,\|\cdot\|_\infty\big).
\end{equation}
By Theorem 3.3.2 of \citet{Edmunds_Triebel_1996}, $\log \cN(\epsilon, A_{R_\kappa}^\infty,\|\cdot\|_\infty)\le C_3 (R_\kappa/\epsilon)^{d/\beta}$, which yields \eqref{eq:SPDE_entropy_bound_app}.
\end{proof}

\begin{lemma}[Centered small-ball probability]\label{lem:SPDE_SBP_centered_app}
Assume $\kappa\ge 1$. There exists $C>0$ (depending only on $d$ and $\beta$) such that for all $\epsilon\in(0,1)$ and all $N\ge 1$,
\begin{equation}\label{eq:SPDE_SBP_centered_app}
    \log P\big(\|f_N\|_\infty\le \epsilon \mid \kappa,N\big)
    \ge
    -C\,\kappa^{d}\,\epsilon^{-d/(\beta-d/2)}.
\end{equation}
\end{lemma}

\begin{proof}
From Proposition 11.11 of \citep{lifshits2012lectures} connecting the metric entropy and the small ball probability, Lemma \ref{lem:SPDE_entropy_H1_app} yields \eqref{eq:SPDE_SBP_centered_app}. 
\end{proof}

We consider $\Besova$ from the main document. For $f^* \in \Besova$ written as
$f^*=\sum_{\bk\in\bbN_0^d} z_{\bk} e_{\bk}$ and $S$ defined in Section 3, we define the spectral cut-off operator at level $J$ as
\begin{equation}\label{eq:SPDE_trunc_operator_app}
    \cT_J f^* = \sum_{\bk\in\bbN_0^d} S(2^{-J}\|\pi\bk\|_2)\,z_{\bk}\,e_{\bk}.
\end{equation}

\begin{lemma}[Sup-norm approximation by $h=\cI_N(\cT_J f^*)$]\label{lem:SPDE_approx_app}
Let $f^*\in\Besova$ and set $C_0=\|f^*\|_{\Besova}$. Let $J$ be the smallest integer such that $2^{-\alpha J}<\epsilon/2$, define $g_J=\cT_J f^*$ and $h=\cI_N g_J$ with $\cT_J$ and $\cI_N$ defined as \eqref{eq:SPDE_trunc_operator_app} and \eqref{eqn: interpolation operator}, respectively. Then $\|f^*-h\|_\infty<\epsilon$ for all sufficiently large $N$. In particular, if $f^*\in L^2(\Omega)$ and
\begin{equation}\label{eq:SPDE_N_condition_app}
    N
    \ge
    C\,\|f^*\|_2^{1/2}\,
    \epsilon^{-\frac12-\frac{2+d/2}{2\alpha}},
\end{equation}
then $\|f^*-h\|_\infty<\epsilon$, for a constant $C>0$ depending only on $(d,\alpha)$.
\end{lemma}

\begin{proof}
By the definition of $\Besova$ and the choice of $J$,
\begin{equation}\label{eq:besov_tail_supp}
    \|f^*-g_J\|_\infty
    =
    \|f^*-\cT_J f^*\|_\infty
    \le
    C_0\,2^{-\alpha J}
    \le
    \epsilon/2.
\end{equation}
It remains to bound $\|g_J-\cI_N g_J\|_\infty$.

Since $g_J=\cT_J f^*$ is a finite cosine series, it is $C^\infty$ on $\Omega$. Therefore, applying
Case 3 in the proof of Lemma 3.2 to $g_J$ yields
\begin{equation}\label{eq:interp_bound_mixed_supp}
    \|g_J-\cI_N g_J\|_\infty
    \le
    \frac{C_1}{N^2}\,
    \max_{|\boldm|=2}\|D^{\boldm} g_J\|_\infty,
\end{equation}
for a constant $C_1>0$ depending only on $d$.

We now bound the second derivatives of $g_J$. From
\begin{equation*}
    g_J = \sum_{\bk\in\mathbb{N}_0^d} a_{\bk} e_{\bk},
    \qquad a_{\bk}=S(2^{-J}\|\pi \bk\|_2)\,z_{\bk},
\end{equation*}
we have $a_{\bk}=0$ whenever $\|\pi \bk\|_2>2^J$ 
since $S$ is supported on $[-1,1]$.
Fix a multi-index $\boldm=(m_1,\ldots,m_d)$ with $|\boldm|=2$. Because
\begin{equation*}
    e_{\bk}(x)=\prod_{r=1}^d e_{k_r}(x_r),
\end{equation*}
and for each $\ell\ge 0$ and $m=0,1,2$,
\begin{equation*}
    \|e_\ell^{(m)}\|_\infty
    \le
    C_2 (\pi \ell)^m
\end{equation*}
holds for some constant $C_2>0$ depending only on $d$, it follows that
\begin{equation*}
    \|D^{\boldm} e_{\bk}\|_\infty
    \le
    C_3 \prod_{r=1}^d (\pi k_r)^{m_r}
    \le
    C_3 \|\pi \bk\|_2^{|\boldm|}
    =
    C_3 \|\pi \bk\|_2^2,
\end{equation*}
for a constant $C_3>0$ depending only on $d$. Therefore, by Cauchy--Schwarz,
\begin{align*}
    \|D^{\boldm} g_J\|_\infty
    &\le
    \sum_{\|\pi \bk\|_2\le 2^J} |a_{\bk}|\,\|D^{\boldm} e_{\bk}\|_\infty \\
    &\le
    C_3
    \Big(
        \sum_{\bk} a_{\bk}^2
    \Big)^{1/2}
    \Big(
        \sum_{\|\pi \bk\|_2\le 2^J} \|\pi \bk\|_2^4
    \Big)^{1/2}.
\end{align*}
Since $0\le S\le 1$,
\begin{equation*}
    \sum_{\bk} a_{\bk}^2
    \le
    \sum_{\bk} z_{\bk}^2
    =
    \|f^*\|_2^2.
\end{equation*}
Also,
\begin{equation*}
    \sum_{\|\pi \bk\|_2\le 2^J} \|\pi \bk\|_2^4
    \le
    2^{4J}\,
    \#\big\{
        \bk\in\mathbb{N}_0^d : \|\pi \bk\|_2\le 2^J
    \big\}
    \le
    C_4\,2^{(4+d)J},
\end{equation*}
for a constant $C_4>0$ depending only on $d$. Hence
\begin{equation*}
    \|D^{\boldm} g_J\|_\infty
    \le
    C_5\,2^{(2+d/2)J}\,\|f^*\|_2,
\end{equation*}
and taking the maximum over $|\boldm|=2$ gives
\begin{equation}\label{eq:derivative_bound_supp}
    \max_{|\boldm|=2}\|D^{\boldm} g_J\|_\infty
    \le
    C_5\,2^{(2+d/2)J}\,\|f^*\|_2.
\end{equation}

Combining \eqref{eq:interp_bound_mixed_supp} and \eqref{eq:derivative_bound_supp},
\begin{equation*}
    \|g_J-\cI_N g_J\|_\infty
    \le
    \frac{C_6}{N^2}\,
    2^{(2+d/2)J}\,\|f^*\|_2.
\end{equation*}
Since $J$ is the smallest integer such that $2^{-\alpha J}<\epsilon/2$, we have
\begin{equation*}
    2^J \le C_7\,\epsilon^{-1/\alpha},
\end{equation*}
for a constant $C_7>0$ depending only on $\alpha$. Therefore,
\begin{equation*}
    \|g_J-\cI_N g_J\|_\infty
    \le
    \frac{C_8\,\|f^*\|_2}{N^2}\,
    \epsilon^{-(2+d/2)/\alpha}.
\end{equation*}
Hence the condition \eqref{eq:SPDE_N_condition_app} ensures that
\begin{equation*}
    \|g_J-\cI_N g_J\|_\infty
    \le
    \epsilon/2.
\end{equation*}
Together with \eqref{eq:besov_tail_supp}, this yields
\begin{equation*}
    \|f^*-h\|_\infty
    \le
    \|f^*-g_J\|_\infty+\|g_J-\cI_N g_J\|_\infty
    <
    \epsilon,
\end{equation*}
which completes the proof.
\end{proof}

\begin{lemma}[RKHS norm bound for the approximant]\label{lem:SPDE_RKHS_norm_app}
Let $f^*\in\Besova$ with $C_0=\|f^*\|_{\Besova}$, and let $h=\cI_N(\cT_J f^*)$ be as in Lemma \ref{lem:SPDE_approx_app}. Then there exist constants $C_1,C_2>0$ (depending only on $d,\alpha,\beta$) such that
\begin{equation}\label{eq:SPDE_RKHS_bound_app}
    \|h\|_{\bbH(\kappa,N)}^2
    \le
    C_1\,\kappa^d
    +
    C_2\,\kappa^{-(2\beta-d)}\,\epsilon^{-(2\beta-2\alpha)/\alpha}.
\end{equation}
\end{lemma}

\begin{proof}
Write
\begin{equation*}
    f^*
    =
    \sum_{\bk\in\mathbb{N}_0^d} z_{\bk} e_{\bk},
    \qquad
    g_J
    =
    \cT_J f^*
    =
    \sum_{\bk\in\mathbb{N}_0^d} a_{\bk} e_{\bk},
    \qquad
    a_{\bk}
    =
    S(2^{-J}\|\pi \bk\|_2)\,z_{\bk}.
\end{equation*}
Since $S$ is supported on $[-1,1]$, we have $a_{\bk}=0$ whenever $\|\pi \bk\|_2>2^J$.
Also, by enlarging the constant in \eqref{eq:SPDE_N_condition_app} if necessary, we may assume
\begin{equation*}
    N\ge 2^J.
\end{equation*}
Hence every index $\bk$ with $\|\pi \bk\|_2\le 2^J$ belongs to $\cJ_N$.

Let $w=\cR_N g_J$. Then
\begin{equation*}
    w
    =
    \sum_{\|\pi \bk\|_2\le 2^J} a_{\bk} v_{\bk}.
\end{equation*}
Since $h=\cI_N g_J=\cE_N w$, Corollary \ref{cor:SPDE_RKHS_ellipsoid_app} yields
\begin{equation}\label{eq:RKHS_start_supp}
    \|h\|_{\bbH(\kappa,N)}^2
    =
    w^\top Q_d w
    =
    \sum_{\|\pi \bk\|_2\le 2^J} \lambda_{\bk}\,\gamma_{\bk}\,a_{\bk}^2.
\end{equation}
Using $\gamma_{\bk}\le 2^d$ and $\lambda_{\bk}\le \Lambda_{\bk}$ from
\eqref{eq:SPDE_lambda_compare_d_app}, we obtain
\begin{equation}\label{eq:RKHS_reduction_supp}
    \|h\|_{\bbH(\kappa,N)}^2
    \le
    2^d
    \sum_{\|\pi \bk\|_2\le 2^J}
    \Lambda_{\bk}\,a_{\bk}^2
    \le
    2^d
    \sum_{\|\pi \bk\|_2\le 2^J}
    \Lambda_{\bk}\,z_{\bk}^2.
\end{equation}

By definition,
\begin{equation*}
    \Lambda_{\bk}
    =
    \kappa^{-(2\beta-d)}
    \big(
        \kappa^2+\|\pi \bk\|_2^2
    \big)^\beta.
\end{equation*}
Using
\begin{equation*}
    (u+v)^\beta \le C_1 (u^\beta+v^\beta),
    \qquad
    u,v\ge 0,
\end{equation*}
we get
\begin{equation*}
    \Lambda_{\bk}
    \le
    C_1 \kappa^{-(2\beta-d)}
    \big(
        \kappa^{2\beta}+\|\pi \bk\|_2^{2\beta}
    \big)
    =
    C_1 \kappa^d
    +
    C_1 \kappa^{-(2\beta-d)}\|\pi \bk\|_2^{2\beta}.
\end{equation*}
Substituting this into \eqref{eq:RKHS_reduction_supp} gives
\begin{equation}\label{eq:RKHS_split_supp}
    \|h\|_{\bbH(\kappa,N)}^2
    \le
    C_2 \kappa^d
    \sum_{\|\pi \bk\|_2\le 2^J} z_{\bk}^2
    +
    C_2 \kappa^{-(2\beta-d)}
    \sum_{\|\pi \bk\|_2\le 2^J} \|\pi \bk\|_2^{2\beta} z_{\bk}^2.
\end{equation}
The first term is bounded by
\begin{equation*}
    C_2 \kappa^d
    \sum_{\bk} z_{\bk}^2
    =
    C_2 \kappa^d \|f^*\|_2^2.
\end{equation*}

For the second term, note that the index $\bk=\mathbf{0}$ contributes $0$. Thus
\begin{equation*}
    \sum_{\|\pi \bk\|_2\le 2^J} \|\pi \bk\|_2^{2\beta} z_{\bk}^2
    =
    \sum_{j=1}^J \sum_{\bk\in A_j} \|\pi \bk\|_2^{2\beta} z_{\bk}^2,
\end{equation*}
where
\begin{equation*}
    A_j
    =
    \big\{
        \bk\in\mathbb{N}_0^d : 2^{j-1}<\|\pi \bk\|_2\le 2^j
    \big\},
    \qquad
    j=1,\ldots,J.
\end{equation*}
Hence
\begin{equation*}
    \sum_{\|\pi \bk\|_2\le 2^J} \|\pi \bk\|_2^{2\beta} z_{\bk}^2
    \le
    \sum_{j=1}^J 2^{2\beta j} \sum_{\bk\in A_j} z_{\bk}^2.
\end{equation*}

Now fix $j\ge 1$. Since
\begin{equation*}
    A_j
    \subset
    \big\{
        \bk\in\mathbb{N}_0^d : \|\pi \bk\|_2>2^{j-1}
    \big\},
\end{equation*}
the definition of $\Besova$ gives
\begin{align*}
    \sum_{\bk\in A_j} z_{\bk}^2
    &\le
    \sum_{\|\pi \bk\|_2>2^{j-1}} z_{\bk}^2 \le
    \|f^*-\cT_{j-1} f^*\|_2^2 \le
    \|f^*-\cT_{j-1} f^*\|_\infty^2 \\
    &\le
    C_0^2\,2^{-2\alpha(j-1)} \le
    C_3\,C_0^2\,2^{-2\alpha j},
\end{align*}
for a constant $C_3>0$ depending only on $\alpha$. Therefore,
\begin{equation*}
    \sum_{\|\pi \bk\|_2\le 2^J} \|\pi \bk\|_2^{2\beta} z_{\bk}^2
    \le
    C_3 C_0^2
    \sum_{j=1}^J 2^{2(\beta-\alpha)j}.
\end{equation*}
Since $\alpha < \beta-d/2<\beta$, we have $\beta-\alpha>0$, and so
\begin{equation*}
    \sum_{j=1}^J 2^{2(\beta-\alpha)j}
    \le
    C_4\,2^{2(\beta-\alpha)J},
\end{equation*}
for a constant $C_4>0$ depending only on $\alpha,\beta$. Substituting into
\eqref{eq:RKHS_split_supp}, we obtain
\begin{equation*}
    \|h\|_{\bbH(\kappa,N)}^2
    \le
    C_5\,\kappa^d
    +
    C_6\,\kappa^{-(2\beta-d)}\,2^{2(\beta-\alpha)J},
\end{equation*}
where $C_5,C_6>0$ absorb $\|f^*\|_2^2$ and $C_0^2$.

Finally, since $J$ is the smallest integer such that $2^{-\alpha J}<\epsilon/2$, we have
\begin{equation*}
    2^J \le C_7\,\epsilon^{-1/\alpha},
\end{equation*}
for a constant $C_7>0$ depending only on $\alpha$. Hence
\begin{equation*}
    2^{2(\beta-\alpha)J}
    \le
    C_8\,\epsilon^{-(2\beta-2\alpha)/\alpha},
\end{equation*}
and therefore
\begin{equation*}
    \|h\|_{\bbH(\kappa,N)}^2
    \le
    C_9\,\kappa^d
    +
    C_{10}\,\kappa^{-(2\beta-d)}\,\epsilon^{-(2\beta-2\alpha)/\alpha}.
\end{equation*}
This proves \eqref{eq:SPDE_RKHS_bound_app}.
\end{proof}

\subsection{Main Proof}\label{appendix:SPDE Main proof}

We verify the three conditions in Theorem 3.3 of \citet{Van2008Rates}. It is enough to construct a sieve $\cB_n \subset C[0, 1]^d$ such that, for $\epsilon_n = n^{-\alpha/(2\alpha+d)}$ up to logarithmic factors,
\begin{align}
    P\big(\|f_N-f^*\|_\infty < \epsilon_n\big)
    &\gtrsim \exp(-\epsilon_n^{-d/\alpha}), \label{eq:spde_prior_conc}\\
    P\big(f_N\notin \cB_n \big)
    &\lesssim \exp(-2\epsilon_n^{-d/\alpha}), \label{eq:spde_sieve}\\
    \log \cN\big(\epsilon_n,\cB_n,\|\cdot\|_n\big)
    &\lesssim \epsilon_n^{-d/\alpha}. \label{eq:spde_entropy}
\end{align}

We call each condition prior concentration, sieve probability, and metric entropy, respectively.

Following the parent GP theory, we assume that the bandwidth parameter $\kappa_n$ satisfies
\begin{equation}\label{eq:SPDE_oracle_seq_app}
    \kappa_n \asymp n^{\frac{\beta-d/2-\alpha}{(2\alpha+d)(\beta-d/2)}}.
\end{equation}
With $\alpha < \beta-d/2$, $\kappa_n \to \infty$ as $n \to \infty$.

\begin{proof}[Proof of the prior concentration condition \eqref{eq:spde_prior_conc}]

First, Let $N_{n, 1}$ be the smallest integer larger than $C\|f^*\|_2^{1/2} \epsilon_n^{-\frac{1}{2} - \frac{2+d/2}{2\alpha}}$. From \eqref{eq:SPDE_N_condition_app}, for every $N \ge N_{n, 1}$, Lemma \ref{lem:SPDE_approx_app} assures that $h = \cI_N (\cT_J f^*)$ satisfies $\|h-f^*\|_\infty \le \epsilon_n$.
Hence, from Lemma \ref{lem:SPDE_SBP_centered_app} and Lemma \ref{lem:SPDE_RKHS_norm_app}, combined with Lemma 5.3 of \cite{van2008reproducing},
\begin{align}
    \log P\big(\|f_N-f^*\|_\infty \le 2\epsilon_n \mid \kappa_n,N\big)
    &\ge \log P(\|f_N\|_{\infty} \le \epsilon_n \mid \kappa_n, N) 
    \\& -\frac{1}{2}\inf_{h: \|h - f^*\|_{\infty} < \epsilon_n}\|h \|_{\bbH(\kappa_n, N)}^2\\
    &\ge
    -C\Big(
        \kappa_n^d \epsilon_n^{-d/(\beta-d/2)}
        + \kappa_n^d
        + \kappa_n^{-(2\beta-d)} \epsilon_n^{-(2\beta-2\alpha)/\alpha}
    \Big)
    \nonumber\\
    &\ge -C' n\epsilon_n^2.
    \label{eq:SPDE_prior_mass_oracle_conditional_app}
\end{align}
The last inequality holds from 
\begin{equation}\label{eq:SPDE_oracle_balance_app}
    \kappa_n^d \epsilon_n^{-d/(\beta-d/2)} \asymp n\epsilon_n^2,
    \qquad
    \kappa_n^{-(2\beta-d)} \epsilon_n^{-(2\beta-2\alpha)/\alpha} \asymp n\epsilon_n^2.
\end{equation}
Therefore, from the lower bound on the prior tail probability of $N$, we obtain
\begin{align}
    P\big(\|f_N-f^*\|_\infty \le 2\epsilon_n\big)
    &\ge
    \inf_{N\ge N_{n, 1}}
    P\big(\|f_N-f^*\|_\infty \le 2\epsilon_n \mid \kappa_n,N\big)\,
    P\big(N\ge N_{n, 1}\big)
    \nonumber\\
    &\ge
    \exp\big(-C_1 n\epsilon_n^2\big).
    \label{eq:SPDE_prior_mass_oracle_app}
\end{align}
This concludes the proof.
\end{proof}

For the proof of \eqref{eq:spde_sieve} and \eqref{eq:spde_entropy}, we first construct $\cB_n$. Let
\begin{equation}\label{eq:SPDE_oracle_RN_app}
    R_n = C_2 \kappa_n^{\beta-d/2}
\end{equation}
and for some $\rho \in (0, \min(\beta - d/2), 1)$, define $N_n$ as the smallest integer satisfying
\begin{equation*}
    N_n \ge \Big(\frac{C_{d, \beta} C_2 \kappa_n^{\beta - d/2} M_n}{\epsilon_n}\Big)^{1/\rho}.
\end{equation*}
From Lemma \ref{lem:holder continuity of AR} and Proof of Lemma 3.2, for all $N \ge N_n$ and $g \in A_{R_n}^{\infty}$,
\begin{equation}\label{eq:g small interpolation error under large N}
    \|g - \cI_N g\|_{\infty} \le \frac{\epsilon_n}{M_n}
\end{equation}
holds. With $M_n > 0$ to be determined, we define
\begin{equation}\label{eq:SPDE_oracle_large_piece_app}
    \cB_n^{(0)}
    :=
    M_n A_{R_n}^\infty + 2\epsilon_n \bbB_1,
\end{equation}
\begin{equation}\label{eq:SPDE_oracle_small_piece_app}
    \cB_{n,N}
    :=
    M_n \bbH_1(\kappa_n,N) + \epsilon_n \bbB_1,
    \qquad N \ge 1,
\end{equation}
and
\begin{equation}\label{eq:SPDE_oracle_sieve_app}
    \cB_n
    :=
    \cB_n^{(0)}
    \cup
    \bigcup_{1\le N < N_n} \cB_{n,N}.
\end{equation}

\begin{proof}[Proof of the sieve condition \eqref{eq:spde_sieve}]
We claim that
\begin{equation}\label{eq:SPDE_oracle_inclusion_app}
    \cB_{n,N} \subset \cB_n,
    \qquad
    \forall\, N\ge 1.
\end{equation}
If $1 \le N < N_n$, this is immediate from \eqref{eq:SPDE_oracle_small_piece_app}. Additionally, for $N \ge N_n$, consider $h, \widetilde h$ defined in the proof of Lemma \ref{lem:SPDE_entropy_H1_app}. As $\widetilde h \in A_{R_n}^{\infty}$ and $\cI_N \widetilde h = h$, from \eqref{eq:g small interpolation error under large N}, we obtain $\|h - \widetilde h\|_\infty \le \epsilon_n / M_n$. This implies that for all $N \ge N_n$ and any $h \in \bbH_1(\kappa_n, N)$, there exists $\widetilde h \in A_{R_n}^\infty$ satisfying $\|h - \widetilde h \|_{\infty} < \epsilon_n / M_n$. As $M_n(h - \widetilde h) \in \epsilon_n \bbB_1$, we obtain
\begin{equation*}
    M_n \bbH_1(\kappa_n, N) + \epsilon_n \bbB_1 \subset \cB_n^{(0)}.
\end{equation*}
Subsequently, for each $N\ge 1$, define
\begin{equation}\label{eq:SPDE_oracle_concentration_fn_app}
    \phi_0^{\kappa_n,N}(\epsilon_n)
    :=
    -\log P\big(\|f_N\|_\infty \le \epsilon_n \mid \kappa_n,N\big).
\end{equation}
By Lemma \ref{lem:SPDE_SBP_centered_app}, there exists $C' > 0$ satisfying
\begin{equation}\label{eq:SPDE_oracle_concentration_bound_app}
    \phi_0^{\kappa_n,N}(\epsilon_n)
    \le
    a_n,
    \qquad
    a_n = C' \kappa_n^d \epsilon_n^{-d/(\beta-d/2)}
    \asymp n\epsilon_n^2.
\end{equation}
Therefore, by Theorem 5.1 of \citet{van2008reproducing},
\begin{align}
    P\big(f_N \notin \cB_{n,N} \mid \kappa_n,N\big)
    &\le
    1-\Phi\Big(
        \Phi^{-1}\big(e^{-\phi_0^{\kappa_n,N}(\epsilon_n)}\big)+M_n
    \Big)
    \nonumber\\
    &\le
    1-\Phi\Big(
        \Phi^{-1}(e^{-a_n})+M_n
    \Big).
    \label{eq:SPDE_oracle_borell_app}
\end{align}
Take $M_n=C_M\sqrt{a_n}$, where $C_M>0$ is a sufficiently large constant. Since $1-\Phi(x)\le \exp(-x^2/2)$ for $x\ge 0$, we have
$\exp(-a_n)\ge 1-\Phi(\sqrt{2a_n})$. Hence
$\Phi^{-1}(\exp(-a_n))\ge -\sqrt{2a_n}$. Taking $C_M\ge 2\sqrt{2}$ gives
$\Phi^{-1}(\exp(-a_n))+M_n\ge M_n/2$. Therefore,
\begin{equation}\label{eq:SPDE_oracle_borell_tail_app}
    P\big(f_N \notin \cB_{n,N} \mid \kappa_n,N\big)
    \le 1- \Phi(M_n / 2) \le
    \exp(-M_n^2/8),
    \qquad
    \forall\, N\ge 1.
\end{equation}
Combined with \eqref{eq:SPDE_oracle_inclusion_app}, we obtain
\begin{equation}\label{eq:SPDE_oracle_sieve_prob_app}
    P(f_N \notin \cB_n)
    \le
    \exp(-M_n^2/8).
\end{equation}
Choosing $C_M$ large enough makes this bound of the required order $\exp(-2n\epsilon_n^2)$.
This concludes the proof.
\end{proof}

Subsequently, we take $M_n=C_M\sqrt{a_n}\asymp n^{\frac{d}{4\alpha + 2d}}$, with $C_M$ fixed sufficiently large.

\begin{proof}[Proof of the metric entropy \eqref{eq:spde_entropy}]
By Theorem 3.3.2 of \citet{Edmunds_Triebel_1996},
\begin{align}
    \log \cN\big(3\epsilon_n,\cB_n^{(0)},\|\cdot\|_\infty\big)
    &\le
    \log \cN\big(\epsilon_n,M_n A_{R_n}^\infty,\|\cdot\|_\infty\big)
    \nonumber\\
    &\le
    C_1 \Big(\frac{M_n R_n}{\epsilon_n}\Big)^{d/\beta}
    \le
    C_2 \kappa_n^{(\beta - d/2)d/\beta}
    \Big(\frac{M_n}{\epsilon_n}\Big)^{d/\beta}.
    \label{eq:SPDE_oracle_entropy_large_app}
\end{align}
Also, for $1 \le N < N_n$, Lemma \ref{lem:SPDE_entropy_H1_app} gives
\begin{align}
    \log \cN\big(2\epsilon_n,\cB_{n,N},\|\cdot\|_\infty\big)
    &\le
    \log \cN\big(\epsilon_n,M_n\bbH_1(\kappa_n,N),\|\cdot\|_\infty\big)
    \nonumber\\
    &=
    \log \cN\big(\epsilon_n/M_n,\bbH_1(\kappa_n,N),\|\cdot\|_\infty\big)\le
    C_3\kappa_n^{(\beta - d/2)d/\beta}
    \Big(\frac{M_n}{\epsilon_n}\Big)^{d/\beta}.
    \label{eq:SPDE_oracle_entropy_small_app}
\end{align}
Hence
\begin{equation}\label{eq:SPDE_oracle_entropy_union_app}
    \log \cN(3\epsilon_n,\cB_n,\|\cdot\|_\infty)
    \le
    C_4 \kappa_n^{(\beta - d/2)d/\beta}
    \Big(\frac{M_n}{\epsilon_n}\Big)^{d/\beta}
    +
    \log(N_n+1).
\end{equation}
Choosing
\begin{equation}\label{eq:SPDE_oracle_Mn_choice_app}
    M_n^2 \asymp n\epsilon_n^2,
\end{equation}
from \eqref{eq:SPDE_oracle_RN_app}, we obtain
\begin{equation}\label{eq:SPDE_oracle_logNn_app}
    \log(N_n+1)
    \le \frac{1}{\rho}\Big\{
        \log M_n + (\beta-d/2)\log \kappa_n + \log(1/\epsilon_n)
    \Big\} + C
    =
    o(n\epsilon_n^2).
\end{equation}
Lastly, plugging in $M_n \asymp n^{\frac{d}{4\alpha + 2d}}$ and $\kappa_n \asymp n^{\frac{\beta - d/2 - \alpha}{(2\alpha + d)(\beta -d/2)}}$, we obtain 
\begin{equation*}
    \log \cN(\epsilon_n, \cB_n, \| \cdot \|_{\infty}) \lesssim n\epsilon_n^2.
\end{equation*}
This concludes the proof.
\end{proof}
\section{Proof of Theorem 3.4}\label{appendix:GPI}
The proof of Theorem 3.4 builds on the parent Gaussian process with squared-exponential covariance kernel \eqref{eqn: RBF kernel function}, whose prior concentration and sieve construction were established in \citet{van2009adaptive}. As introduced in Section \ref{ssec:GPI}, GPI method can be represented as $f_N=\cI_N g$ where $g$ is a parent GP with mean zero and squared-exponential covariance kernel, and $\cI_N$ is an interpolation operator defined in \eqref{eqn: interpolation operator}. Thus, the central goal is to control the interpolation error $\|g-\cI_N g\|_\infty$ with high probability under the law of parent GP. Section \ref{appendix:GPI_technical} collects the technical lemmas that bound this error via modulus of continuity of parent GP, and Sections \ref{appendix:GPI_main} and \ref{appendix:gpi_sieve_universalN} combine these with the results of \citet{van2009adaptive} to verify the three contraction conditions, as in Section \ref{appendix:SPDE Main proof}.
\subsection{Technical lemmas for the GPI approach}\label{appendix:GPI_technical}

Let $g=\{g(s):s\in[0,1]^d\}$ be a mean-zero squared-exponential GP with covariance
\begin{equation}\label{eqn:gpi_covariance}
    \bbE\{g(s)g(t)\}=\exp\big\{-\kappa^2\|s-t\|_2^2\big\},
\end{equation}
where $\kappa>0$ is the inverse bandwidth parameter.
Define the canonical pseudometric
\begin{equation}\label{eq:dg_def_supp}
    d_g(s,t)=\Big[\bbE\{g(s)-g(t)\}^2\Big]^{1/2}.
\end{equation}
For $\epsilon>0$, define
\begin{align}
B^\infty_\epsilon&=\Big\{(s,t)\in[0,1]^d\times[0,1]^d:\|s-t\|_\infty\le \epsilon\Big\},\nonumber
\\
B^g_\epsilon&=\Big\{(s,t)\in[0,1]^d\times[0,1]^d:d_g(s,t)\le \epsilon\Big\}.\label{eq:B_sets_supp}
\end{align}

Let $\{u_{j_r}=j_r/N:0\le j_r\le N\}$ be the regular grid and $\{\psi_{j_r}\}$ be the hat basis in each coordinate.
Define the $d$-dimensional interpolation operator by
\begin{equation}\label{eq:gpi_interp_operator_aux}
    (\cI_N h)(x_1,\ldots,x_d)
    =\sum_{j\in\cJ_N} h(u_{j_1},\ldots,u_{j_d})\,
    \psi_{j_1}(x_1)\cdots\psi_{j_d}(x_d),
\end{equation}
and set $f_N=\cI_N g$.

\begin{lemma}\label{lem:gpi_metric_relation}
For all $s,t\in[0,1]^d$,
\begin{equation}\label{eq:gpi_metric_relation_eq}
    d_g(s,t)^2
    =2\Big\{1-\exp\big(-\kappa^2\|s-t\|_2^2\big)\Big\}
    \le 2\kappa^2\|s-t\|_2^2
    \le 2\kappa^2 d \|s-t\|_\infty^2.
\end{equation}
In particular, $B^\infty_\epsilon \subset B^g_{\kappa\sqrt{2d}\,\epsilon}$.
\end{lemma}
\begin{proof}
    The result follows from $1-\exp(-x) \le x$ for $x \ge 0$ and $\|v\|_2^2 \le d\|v\|_\infty^2$ for $v \in \bbR^d$.
\end{proof}

\begin{lemma}\label{lem:gpi_fluctuation_expectation}
There exists a universal constant $K>0$ such that for all $\kappa>0$ and all integers $N\ge 2$,
\begin{equation}\label{eq:gpi_Esup_modulus}
    \bbE \sup_{(s,t)\in B^\infty_{1/N}} \{g(s)-g(t)\}
    \le \frac{C_{\mathrm{osc}}\,\kappa d}{N}
    \Big(\log(C_{\mathrm{vol}} N)\Big)^{1/2},
\end{equation}
where $C_{\mathrm{osc}},C_{\mathrm{vol}}>0$ are constants depending only on $d$.
\end{lemma}

\begin{proof}
By Lemma \ref{lem:gpi_metric_relation},
\begin{equation}\label{eq:gpi_Esup_step1}
    \bbE\sup_{(s,t)\in B^\infty_{1/N}} \{g(s)-g(t)\}
    \le \bbE\sup_{(s,t)\in B^g_{\kappa\sqrt{2d}/N}} \{g(s)-g(t)\}.
\end{equation}
Following the argument used in Theorem 4.5 of \citet{adler1990introduction}, there exists a universal constant $K>0$ such that
\begin{align}
    &\bbE\sup_{(s,t)\in B^g_{\kappa\sqrt{2d}/N}} \{g(s)-g(t)\} \nonumber
    \\&\le 8K \sup_{t\in[0,1]^d} \int_0^{\kappa\sqrt{2d}/(2N)}
    \Bigg[\log\frac{1}{\mbox{Vol}\big\{B(t,\epsilon,\|\cdot\|_g)\cap[0,1]^d\big\}}\Bigg]^{1/2} d\epsilon. \label{eq:gpi_chaining_integral}
\end{align}
Using Lemma \ref{lem:gpi_metric_relation}, $B(t,\epsilon/(\kappa\sqrt{2}),\|\cdot\|_2)\cap[0,1]^d \subset B(t,\epsilon,\|\cdot\|_g)\cap[0,1]^d$.
Moreover, there exists a constant $c_d\in(0,1]$ such that for all $t\in[0,1]^d$ and all $r\in(0,1]$,
\begin{equation}\label{eq:gpi_vol_lower}
    \mbox{Vol}\Big\{B\big(t,r,\|\cdot\|_2\big)\cap[0,1]^d\Big\}\ge c_d v_d r^d,
\end{equation}
where $v_d$ is the volume of the unit Euclidean ball in $\bbR^d$.
Substituting these bounds into the integral and evaluating as in the standard chaining calculation yields \eqref{eq:gpi_Esup_modulus} with constants depending only on $d$.
\end{proof}

\begin{lemma}\label{lem:gpi_fluctuation_tail}
Let $N\ge 2$. If $\kappa$ and $N$ satisfy
\begin{equation}\label{eq:gpi_modulus_condition}
    \frac{C_{\mathrm{osc}}\,\kappa d}{N}
    \Big(\log(C_{\mathrm{vol}} N)\Big)^{1/2} \le \epsilon/4,
\end{equation}
then
\begin{equation}\label{eq:gpi_modulus_tail}
    P\Bigg(\sup_{(s,t)\in B^\infty_{1/N}} |g(s)-g(t)| > \epsilon/2 \,\Big|\, \kappa\Bigg)
    \le 4\exp\Big(-c_{\mathrm{Bor}}\frac{N^2\epsilon^2}{\kappa^2 d}\Big),
\end{equation}
for a universal constant $c_{\mathrm{Bor}}>0$.
\end{lemma}

\begin{proof}
Define the centered Gaussian process $Y_{(s,t)}=g(s)-g(t)$ indexed by $B^\infty_{1/N}$.
By Lemma \ref{lem:gpi_metric_relation},
\begin{equation}\label{eq:gpi_var_bound}
    \sup_{(s,t)\in B^\infty_{1/N}} \mbox{Var}\big\{Y_{(s,t)}\big\}
    =\sup_{(s,t)\in B^\infty_{1/N}} d_g(s,t)^2
    \le \frac{2\kappa^2 d}{N^2}.
\end{equation}
By Lemma \ref{lem:gpi_fluctuation_expectation} and \eqref{eq:gpi_modulus_condition},
\begin{equation}\label{eq:gpi_Esup_bound}
    \bbE\sup_{(s,t)\in B^\infty_{1/N}} Y_{(s,t)} \le \epsilon/4.
\end{equation}
Borell's inequality \citep[Theorem 2.1]{adler1990introduction} implies that for any $u>0$,
\begin{equation}\label{eq:gpi_borell}
    P\Bigg(\sup_{(s,t)\in B^\infty_{1/N}} Y_{(s,t)} > \bbE\sup_{(s,t)\in B^\infty_{1/N}} Y_{(s,t)} + u \,\Big|\,\kappa\Bigg)
    \le \exp\Big(-\frac{u^2}{2\sigma^2}\Big),
\end{equation}
where $\sigma^2 \le 2\kappa^2 d/N^2$ by \eqref{eq:gpi_var_bound}. Taking $u=\epsilon/4$ yields
\begin{equation}\label{eq:gpi_one_sided_tail}
    P\Bigg(\sup_{(s,t)\in B^\infty_{1/N}} Y_{(s,t)} > \epsilon/2 \,\Big|\,\kappa\Bigg)
    \le \exp\Big(-c_{\mathrm{Bor}}\frac{N^2\epsilon^2}{\kappa^2 d}\Big).
\end{equation}
Using symmetry of $Y_{(s,t)}$ and a union bound gives \eqref{eq:gpi_modulus_tail}.
\end{proof}

\subsection{Main proof}\label{appendix:GPI_main}

Based on Theorem 3.3 of \citet{Van2008Rates}, it suffices to verify the following three conditions. There exists a sieve $\cF^*(\epsilon)\subset C[0,1]^d$ and a constant $C>0$ such that for all sufficiently small $\epsilon>0$,
\begin{align}
    P\big(\|f_N-f^*\|_\infty < \epsilon\big)
    &\gtrsim \exp\Big\{-\epsilon^{-d/\alpha}\log^C(1/\epsilon)\Big\}, \label{eq:gpi_prior_conc}\\
    P\big(f_N\notin \cF^*(\epsilon)\big)
    &\lesssim \exp\Big\{-2\epsilon^{-d/\alpha}\log^C(1/\epsilon)\Big\}, \label{eq:gpi_sieve}\\
    \log \cN\big(\epsilon,\cF^*(\epsilon),\|\cdot\|_n\big)
    &\lesssim \epsilon^{-d/\alpha}\log^C(1/\epsilon). \label{eq:gpi_entropy}
\end{align}

\begin{proof}[Proof of the prior concentration condition \eqref{eq:gpi_prior_conc}]
Fix $\epsilon\in(0,\epsilon_0)$ and let $f_N=\cI_N g$ with $g$ defined in \eqref{eqn:gpi_covariance}.
We lower bound $P(\|f_N-f^*\|_\infty<\epsilon)$ by restricting to suitable $(\kappa,N)$.

Let $\kappa_\epsilon=(C_0/\epsilon)^{1/\alpha}$, where $C_0>0$ is the constant appearing in the small ball condition of \citet{van2009adaptive}.
For $\kappa\in[\kappa_\epsilon,2\kappa_\epsilon]$ and any integer $N\ge 2$, by the triangle inequality,
\begin{align}
    P\big(\|f_N-f^*\|_\infty<\epsilon \,\big|\, \kappa,N\big)
    &\ge P\big(\|g-f^*\|_\infty<\epsilon/2 \,\big|\,\kappa\big)
    \\&- P\big(\|g-\cI_N g\|_\infty>\epsilon/2 \,\big|\,\kappa,N\big). \label{eq:gpi_triangle}
\end{align}

For the first term, Section 5.1 of \citet{van2009adaptive} implies that for all $\kappa\ge \kappa_0$ and all $\epsilon$ satisfying $\epsilon>C_0\kappa^{-\alpha}$,
\begin{equation}\label{eq:gpi_smallball}
    P\big(\|g-f^*\|_\infty<\epsilon/2 \,\big|\,\kappa\big)
    \ge \exp\Big[-K_1 \kappa^d \Big\{\log\Big(\frac{2\kappa}{\epsilon}\Big)\Big\}^{1+d}\Big],
\end{equation}
for constants $K_1,\kappa_0>0$.

For the second term, Lemma \ref{lem:gpi_fluctuation_tail} implies
\begin{equation}\label{eq:gpi_interp_error_control}
    \|g-\cI_N g\|_\infty
    \le \sup_{(s,t)\in B^\infty_{1/N}} |g(s)-g(t)|.
\end{equation}
Choose $N_1(\epsilon)$ as any integer satisfying, with $\kappa=2\kappa_\epsilon$,
\begin{equation}\label{eq:gpi_N1_choice}
    \frac{C_{\mathrm{osc}}\,(2\kappa_\epsilon) d}{N_1(\epsilon)}
    \Big(\log\big(C_{\mathrm{vol}} N_1(\epsilon)\big)\Big)^{1/2} \le \epsilon/4,
\end{equation}
and
\begin{equation}\label{eq:gpi_N1_exponent}
    c_{\mathrm{Bor}}\frac{N_1(\epsilon)^2\epsilon^2}{(2\kappa_\epsilon)^2 d}
    \ge 2\epsilon^{-d/\alpha}\log^C(1/\epsilon).
\end{equation}
Then for all $\kappa\in[\kappa_\epsilon,2\kappa_\epsilon]$ and all $N\ge N_1(\epsilon)$,
\begin{equation}\label{eq:gpi_interp_prob_control}
    P\big(\|g-\cI_N g\|_\infty>\epsilon/2 \,\big|\,\kappa,N\big)
    \le 4\exp\Big\{-2\epsilon^{-d/\alpha}\log^C(1/\epsilon)\Big\}.
\end{equation}

Combining \eqref{eq:gpi_triangle}, \eqref{eq:gpi_smallball}, and \eqref{eq:gpi_interp_prob_control}, we obtain for all sufficiently small $\epsilon$,
\begin{equation}\label{eq:gpi_conditional_conc_lower}
    P\big(\|f_N-f^*\|_\infty<\epsilon \,\big|\,\kappa,N\big)
    \ge \frac{1}{2}\exp\Big[-K_1 \kappa^d \Big\{\log\Big(\frac{2\kappa}{\epsilon}\Big)\Big\}^{1+d}\Big]
\end{equation}
uniformly over $\kappa\in[\kappa_\epsilon,2\kappa_\epsilon]$ and $N\ge N_1(\epsilon)$.

Integrating over $(\kappa,N)$ and using the prior distribution on $\kappa$ and $N$ yields
\begin{align}\label{eq:gpi_prior_mass_final}
    P\big(\|f_N-f^*\|_\infty<\epsilon\big)
    &\ge P_2\big(N\ge N_1(\epsilon)\big)
    \int_{\kappa_\epsilon}^{2\kappa_\epsilon}
    \frac{1}{2}\exp\Big[-K_1 \kappa^d \Big\{\log\Big(\frac{2\kappa}{\epsilon}\Big)\Big\}^{1+d}\Big] P_1(d\kappa) \\
    &\gtrsim \exp\Big\{-\epsilon^{-d/\alpha}\log^C(1/\epsilon)\Big\},
\end{align}
after adjusting constants. This proves \eqref{eq:gpi_prior_conc}.
\end{proof}

\subsection{Sieve construction for the GPI prior}\label{appendix:gpi_sieve_universalN}

We construct $\cF^*(\epsilon)$ in three steps:
\begin{enumerate}
    \item choose cutoffs $\kappa_\epsilon$ and $N_0(\epsilon)$ so that a modulus-of-continuity event at scale $1/N_0(\epsilon)$ has exponentially large probability uniformly over $\kappa\le \kappa_\epsilon$;
    \item stabilize the sieve $\cG(\epsilon)$ of \citet{van2009adaptive} by intersecting it with this modulus event;
    \item define $\cF^*(\epsilon)$ so that it covers $f_N$ for all $N\ge 1$ without requiring an upper tail bound on $P_2(N>m)$.
\end{enumerate}

The following basic properties of $\cI_N$ will be used repeatedly.

\begin{lemma}\label{lem:gpi_interp_basic}
For any $h_1,h_2\in C[0,1]^d$ and any integer $N\ge 1$,
\begin{align}
    \|\cI_N h_1 - \cI_N h_2\|_\infty
    &\le \|h_1-h_2\|_\infty, \label{eq:gpi_interp_contraction}\\
    \|h-\cI_N h\|_\infty
    &\le \sup_{(s,t)\in B^\infty_{1/N}} |h(s)-h(t)|. \label{eq:gpi_interp_error_modulus}
\end{align}
\end{lemma}

\begin{proof}
Fix $x\in[0,1]^d$. The weights $\psi_{j_1}(x_1)\cdots\psi_{j_d}(x_d)$ are nonnegative and sum to one, so
\begin{equation}\label{eq:gpi_contraction_pf}
    |(\cI_N h_1)(x)-(\cI_N h_2)(x)|
    \le \|h_1-h_2\|_\infty \sum_{j\in\cJ_N}\psi_{j_1}(x_1)\cdots\psi_{j_d}(x_d)
    =\|h_1-h_2\|_\infty,
\end{equation}
which proves \eqref{eq:gpi_interp_contraction}.
For \eqref{eq:gpi_interp_error_modulus}, $(\cI_N h)(x)$ is a convex combination of values $h(u)$ at vertices $u$ with $\|x-u\|_\infty\le 1/N$, hence
\begin{equation}\label{eq:gpi_modulus_pf}
    |h(x)-(\cI_N h)(x)|
    \le \sup_{(s,t)\in B^\infty_{1/N}} |h(s)-h(t)|.
\end{equation}
Taking the supremum over $x$ yields \eqref{eq:gpi_interp_error_modulus}.
\end{proof}

Choose $\kappa_\epsilon$ so that
\begin{equation}\label{eq:gpi_kappa_cutoff}
    P_1(\kappa>\kappa_\epsilon)
    \le \exp\Big\{-2\epsilon^{-d/\alpha}\log^C(1/\epsilon)\Big\}.
\end{equation}
Let $N_0(\epsilon)$ be any integer satisfying, with $\kappa=\kappa_\epsilon$,
\begin{equation}\label{eq:gpi_N0_choice}
    \frac{C_{\mathrm{osc}}\,\kappa_\epsilon d}{N_0(\epsilon)}
    \Big(\log\big(C_{\mathrm{vol}} N_0(\epsilon)\big)\Big)^{1/2}\le \epsilon/4,
\end{equation}
and
\begin{equation}\label{eq:gpi_N0_exponent}
    c_{\mathrm{Bor}}\frac{N_0(\epsilon)^2\epsilon^2}{\kappa_\epsilon^2 d}
    \ge 2\epsilon^{-d/\alpha}\log^C(1/\epsilon).
\end{equation}

Define the modulus event
\begin{equation}\label{eq:gpi_modulus_event}
    \cE(\epsilon)=\Big\{g:\sup_{(s,t)\in B^\infty_{1/N_0(\epsilon)}} |g(s)-g(t)| \le \epsilon/2\Big\}.
\end{equation}
By Lemma \ref{lem:gpi_fluctuation_tail}, \eqref{eq:gpi_N0_choice}, and \eqref{eq:gpi_N0_exponent}, for all $\kappa\le \kappa_\epsilon$,
\begin{equation}\label{eq:gpi_modulus_event_prob}
    P\big(g\notin \cE(\epsilon)\,\big|\,\kappa\big)
    \le 4\exp\Big\{-2\epsilon^{-d/\alpha}\log^C(1/\epsilon)\Big\}.
\end{equation}

Let $\cG(\epsilon)$ be the sieve defined in Section 5 of \citet{van2009adaptive} with parameters $(M,r,\delta)$ chosen so that
\begin{equation}\label{eq:gpi_baseG_prob}
    P\big(g\notin \cG(\epsilon)\big)
    \le \exp\Big\{-2\epsilon^{-d/\alpha}\log^C(1/\epsilon)\Big\}.
\end{equation}
Define the stabilized sieve
\begin{equation}\label{eq:gpi_stabilizedG_def}
    \widetilde{\cG}(\epsilon)=\cG(\epsilon)\cap \cE(\epsilon).
\end{equation}

Finally, define the universal-$N$ sieve for $\{f_N\}_{N\ge 1}$ by
\begin{equation}\label{eq:gpi_universalF_def}
    \cF^*(\epsilon)
    = \big(\widetilde{\cG}(\epsilon)+\epsilon \bbB_1\big)
    \cup \bigcup_{1\le N\le N_0(\epsilon)} \cI_N\big(\widetilde{\cG}(\epsilon)\big),
\end{equation}
where $\bbB_1$ denotes the unit ball of $C[0,1]^d$ under $\|\cdot\|_\infty$.

\begin{proof}[Proof of the sieve condition \eqref{eq:gpi_sieve}]
Let $f_N=\cI_N g$.
If $g\in\widetilde{\cG}(\epsilon)$ and $N\le N_0(\epsilon)$, then $f_N\in \cI_N(\widetilde{\cG}(\epsilon))\subset \cF^*(\epsilon)$ by \eqref{eq:gpi_universalF_def}.
If $g\in\widetilde{\cG}(\epsilon)$ and $N\ge N_0(\epsilon)$, then $g\in\cE(\epsilon)$ and by \eqref{eq:gpi_interp_error_modulus},
\begin{equation}\label{eq:gpi_sieve_case2}
    \|f_N-g\|_\infty
    \le \sup_{(s,t)\in B^\infty_{1/N_0(\epsilon)}} |g(s)-g(t)|
    \le \epsilon/2,
\end{equation}
so $f_N\in \widetilde{\cG}(\epsilon)+\epsilon\bbB_1\subset \cF^*(\epsilon)$.
Therefore, for all $N\ge 1$,
\begin{equation}\label{eq:gpi_sieve_inclusion}
    \{g\in\widetilde{\cG}(\epsilon)\}\subset \{f_N\in \cF^*(\epsilon)\}.
\end{equation}
Hence
\begin{equation}\label{eq:gpi_sieve_bound_start}
    P\big(f_N\notin \cF^*(\epsilon)\big)\le P\big(g\notin \widetilde{\cG}(\epsilon)\big).
\end{equation}
Using \eqref{eq:gpi_stabilizedG_def} and a union bound,
\begin{equation}\label{eq:gpi_union_bound}
    P\big(g\notin \widetilde{\cG}(\epsilon)\big)
    \le P\big(g\notin \cG(\epsilon)\big) + P\big(g\notin \cE(\epsilon)\big).
\end{equation}
The first term is controlled by \eqref{eq:gpi_baseG_prob}.
For the second term,
\begin{equation}\label{eq:gpi_E_bound}
    P\big(g\notin \cE(\epsilon)\big)
    \le P_1(\kappa>\kappa_\epsilon)+\sup_{\kappa\le \kappa_\epsilon}P\big(g\notin \cE(\epsilon)\,\big|\,\kappa\big),
\end{equation}
which is bounded by \eqref{eq:gpi_kappa_cutoff} and \eqref{eq:gpi_modulus_event_prob}.
Combining these bounds yields \eqref{eq:gpi_sieve} after adjusting constants.
\end{proof}

\begin{proof}[Proof of the metric entropy condition \eqref{eq:gpi_entropy}]
Since $\|h\|_n\le \|h\|_\infty$, it suffices to work under $\|\cdot\|_\infty$.
We first use the shift bound
\begin{equation}\label{eq:gpi_entropy_shift}
    \cN\big(2\epsilon,\widetilde{\cG}(\epsilon)+\epsilon\bbB_1,\|\cdot\|_\infty\big)
    \le \cN\big(\epsilon,\widetilde{\cG}(\epsilon),\|\cdot\|_\infty\big)
    \le \cN\big(\epsilon,\cG(\epsilon),\|\cdot\|_\infty\big).
\end{equation}
Next, by \eqref{eq:gpi_interp_contraction},
\begin{equation}\label{eq:gpi_entropy_contraction}
    \cN\big(\epsilon,\cI_N(\widetilde{\cG}(\epsilon)),\|\cdot\|_\infty\big)
    \le \cN\big(\epsilon,\widetilde{\cG}(\epsilon),\|\cdot\|_\infty\big)
    \le \cN\big(\epsilon,\cG(\epsilon),\|\cdot\|_\infty\big).
\end{equation}
Using \eqref{eq:gpi_universalF_def} and a union bound for covering numbers,
\begin{align}\label{eq:gpi_entropy_union}
    \cN\big(2\epsilon,\cF^*(\epsilon),\|\cdot\|_\infty\big)
    &\le \cN\big(2\epsilon,\widetilde{\cG}(\epsilon)+\epsilon\bbB_1,\|\cdot\|_\infty\big)
    +\sum_{N=1}^{N_0(\epsilon)} \cN\big(2\epsilon,\cI_N(\widetilde{\cG}(\epsilon)),\|\cdot\|_\infty\big) \\
    &\le \big(1+N_0(\epsilon)\big)\cN\big(\epsilon,\cG(\epsilon),\|\cdot\|_\infty\big).
\end{align}
Taking logarithms yields
\begin{equation}\label{eq:gpi_entropy_log}
    \log \cN\big(2\epsilon,\cF^*(\epsilon),\|\cdot\|_\infty\big)
    \le \log\big(1+N_0(\epsilon)\big)+\log \cN\big(\epsilon,\cG(\epsilon),\|\cdot\|_\infty\big).
\end{equation}
From Section 5 of \citet{van2009adaptive}, $\log \cN(\epsilon,\cG(\epsilon),\|\cdot\|_\infty)\lesssim \epsilon^{-d/\alpha}\log^C(1/\epsilon)$.
Moreover, \eqref{eq:gpi_N0_choice}--\eqref{eq:gpi_N0_exponent} imply that $N_0(\epsilon)$ grows at most polynomially in $\kappa_\epsilon$ and $\epsilon^{-1}$ up to a logarithmic factor, hence $\log(1+N_0(\epsilon))\lesssim \log(1/\epsilon)$.
This proves \eqref{eq:gpi_entropy}.
\end{proof}

\bibliography{reference_arxiv_260606}
\bibliographystyle{plainnat}
\end{document}